\documentclass[a4paper,reqno,12pt]{amsart}

\usepackage{rotating,pdflscape,color,amssymb,hyperref,subfigure,psfrag,amsmath,eufrak,bbm,dsfont}

\title[]{Large deviations of
the  extreme eigenvalues of 
random deformations of  matrices}



\date{\today}

\newcommand{\Hr}{\mathsf{H}_r}

\newcommand{\ovl}{\overline}

\newcommand{\bbm}{\begin{bmatrix}}

\newcommand{\ebm}{\end{bmatrix}}

\newcommand{\bes}{\begin{equation*}}

\newcommand{\ees}{\end{equation*}}

\newcommand{\be}{\begin{equation}}

\newcommand{\ee}{\end{equation}}

\newcommand{\beqy}{\begin{eqnarray}}

\newcommand{\eeqy}{\end{eqnarray}}

\newcommand{\beq}{\begin{eqnarray*}}

\newcommand{\eeq}{\end{eqnarray*}}

\newcommand{\lan}{\langle}

\newcommand{\ran}{\rangle}

\newcommand{\diag}{\operatorname{diag}}

\newcommand{\e}{\varepsilon}
\newcommand{\Pro}{\mathbb{P}}

\newcommand{\tr}{\operatorname{tr}}

\newcommand{\Tr}{\operatorname{Tr}}

\newcommand{\ninf}{\underset{n\to+\infty}{\longrightarrow}}

\newcommand{\one}{\mathbbm{1}}

\newcommand{\E}{\mathbb{E}}

\newcommand{\R}{\mathbb{R}}

\newcommand{\C}{\mathbb{C}}

\newcommand{\n}{\mathbb{N}}

\newcommand{\K}{\mathcal{K}}

\newcommand{\ud}{\mathrm{d}}

\newcommand{\pro}{probability }

\newcommand{\f}{\frac}

\newcommand{\ff}{\frac{1}}

\newcommand{\lf}{\left}

\newcommand{\ri}{\right}

\newcommand{\st}{such that }

\newcommand{\la}{\lambda}

\newcommand{\La}{\Lambda}

\newcommand{\ste}{\, ;\, }

\newcommand{\mc}{\mathcal }

\newcommand{\eps}{\varepsilon}

\def\ra{{\rightarrow}}
\newcommand{\wtX}{\widetilde{X_n}}

\newcommand{\wtl}{\widetilde{\la}}

\newcommand{\bck}{\backslash}

\newtheorem{Th}{Theorem}[section]

\newtheorem{assum}[Th]{Assumption}
\newtheorem{propo}[Th]{Proposition}

\newtheorem{propr}[Th]{Property}

\newtheorem{lem}[Th]{Lemma}

\newtheorem{rmq}[Th]{Remark}

\newtheorem{cor}[Th]{Corollary}

\newenvironment{pr}{\noindent {\it Proof. }}{\hfill$\square$}

\long\def\symbolfootnote[#1]#2{\begingroup

\def\thefootnote{\fnsymbol{footnote}}\footnote[#1]{#2}\endgroup}
\author[F. Benaych-Georges, A. Guionnet, M. Maida]{F. Benaych-Georges*, A. Guionnet$^\star$, M. Maida$^\sharp$.}
\thanks {* UPMC Univ Paris 6, LPMA,  Case courier 188, 4, Place Jussieu, 75252 Paris Cedex 05, France, florent.benaych@upmc.fr,\\
$\star$ UMPA, ENS Lyon, 46 all\'ee d'Italie, 69364 Lyon Cedex 07,
France, aguionne@umpa.ens-lyon.fr,\\
$\sharp$ Universit\'e Paris-Sud, Laboratoire de Math\'ematiques,
B\^atiment 425,
Facult\'e des Sciences,
91405 Orsay Cedex, France, mylene.maida@math.u-psud.fr.\\
This work was   supported by the \emph{Agence Nationale de la
Recherche} grant ANR-08-BLAN-0311-03.}

\begin{document}

\begin{abstract}
 Consider a real diagonal deterministic matrix $X_n$ of size $n$ with spectral measure converging to a compactly supported probability measure. We perturb this matrix by adding a random finite rank matrix, with delocalized eigenvectors. We show that the joint law of the extreme eigenvalues of the perturbed model satisfies a large deviation principle  in the scale $n$, with a good rate function given by
a variational formula. \\
We tackle both cases when the extreme eigenvalues of $X_n$ converge to the edges of the support of the limiting measure and when we allow some eigenvalues of $X_n$, that we call outliers, to converge out of the bulk.\\
We can also generalise our results to the case when $X_n$ is random, with law proportional to $e^{- n\Tr V(X)}\ud X,$ for  $V$ growing fast enough at infinity and any perturbation of finite rank.
\end{abstract}

\maketitle

\vspace*{-0.5cm}
\tableofcontents
\section{Introduction}\label{SectionDefmodel}
In the last twenty years, many features of
the asymptotics of the spectrum of large random matrices
have been understood. For a wide variety of classical
models of random matrices (the canonical examples hereafter will
be 
 Wigner matrices  \cite{Wigner}, or    Wishart matrices \cite{pastur-marchenko}),
it has been shown that the spectral measure converges 
almost surely.
The extreme eigenvalues converge for most of these models to the boundaries
of the limiting spectral measure (see e.g. \cite{furedi} or \cite{bai2}).
Fluctuations  of  the
spectral measure and the extreme eigenvalues of these
models  could also 
be studied under a fair generality
over the entries of the matrices; we refer to  \cite{soshnikov}
and \cite{greg-ofer}, or  \cite{alice-greg-ofer} and
\cite{Bai-Silverstein} for reviews. Recently, even the fluctuations
of the eigenvalues inside the bulk could be studied for
rather general   entries  and were shown to be universal (see e.g. \cite{eprsy}
or \cite{TaoVu}). Concentration of measure phenomenon
and moderate deviations could also be established in \cite{GZconc, DGZmod,mod}.\\

 Yet,
the understanding of the large deviations
of the spectrum of large random matrices is still very scarce
and exists only in very specific 
cases. Indeed, the spectrum of a matrix is a very complicated 
function of the entries, so that usual large deviation
theorems, mainly based on independence, do not apply.
Moreover, large deviations rate functions
have to depend on the distribution 
of the entries and only guessing 
their definition is still  a widely open
question.  In the case of Gaussian Wigner matrices, where the joint law of the eigenvalues is simply
given by a Coulomb gas Gibbs measure, 
things are much easier  and a full large deviation
principle for the law of the spectral measure of such
matrices 
was proved in  \cite{BAG3}. This extends
to other ensembles  distributed according to similar
Gibbs measure, for instance  Gaussian Wishart
matrices   \cite{hiaipetz}. Similar large deviation
results hold in 
 discrete situations with a 
Coulomb gas distribution \cite{delph}. A large deviation
principle was also established in \cite{GZ}
for the law of the spectral measure of 
a random matrix given as 
  the sum of a self-adjoint Gaussian  Wigner random matrix
and a deterministic self-adjoint matrix (or as a Gaussian
Wishart matrix with non trivial
covariance matrix). In this case, the proof 
uses stochastic analysis and Dyson's Brownian motion, as there is
no explicit joint law for the eigenvalues, but again relies
heavily on the fact that the random matrix has
Gaussian entries.\\
The large deviations for the law of the extreme eigenvalues
were studied in a slightly more general setting. Again relying on
the explicit joint law of the eigenvalues, a large deviation 
principle was derived in \cite{bdg01} for the same Gaussian type
 models.   
   The large deviations of
extreme eigenvalues of Gaussian Wishart matrices were studied
in \cite{Boh}.  
In the case where the Wishart
matrix is of the form $XX^*$ with $X$ a $n\times r$ rectangular
matrix so that the ratio $r/n$ of its dimensions 
goes to zero, large deviations bounds for the extreme
eigenvalues could be derived under
more general
assumptions on the entries in \cite{Fey}.
Our approaches allow also to obtain  a
full large deviation for the spectrum of such
Wishart matrices when $r$ is kept fixed while $n$ goes
to infinity (see Section \ref{wishartsec}).\\

 In this
article, we shall be concerned with the effect of finite rank
deformations on the deviations of the extreme eigenvalues
of random matrices. In fact, using Weyl's interlacing property, it is easy to check that such
finite rank perturbations
do not change the deviations
of the spectral measure. 
But it strongly affects the behavior of a few extreme eigenvalues,
not only at the level of deviations but also as far as convergence and fluctuations are concerned. In the case
of  Gaussian Wishart matrices, the asymptotics of these extreme
eigenvalues were
established in \cite{BBP}
 and a sharp phase transition, known
as the BBP transition, was exhibited. According to
the strength of the perturbation, the extreme eigenvalues
converge  to the edge of the bulk or away from the bulk.
The fluctuations of these eigenvalues
were also shown in \cite{BBP} to be  given either  
by the Tracy-Widom distribution in the first
case, or by the Gaussian distribution in the
second case. Universality (and non-universality)
 of the fluctuations in BBP transition  was studied for various models,
see e.g. \cite{CDF09,CDF09b,FP07,bai-yao-TCL}.

The goal of this article is to study 
the large deviations of the extreme
eigenvalues of  such finite rank perturbations
of large random matrices. 
 In \cite{myleneEJP},
a large deviation principle for the largest 
eigenvalue of matrices of the GOE and GUE
 deformed by a rank one matrix was obtained 
by using fine asymptotics of the Itzykson-Zuber-Harich-Chandra 
(or spherical) integrals.
The large deviations
of the extreme eigenvalues of
a Wigner matrix perturbed by a matrix with finite rank  greater
than one
happened to be much more complicated. One of the outcomes
of this paper is to prove such a large deviation result
when the Wigner matrix is Gaussian. In fact, our 
result will include the more general case 
where the non-perturbed matrices are taken in some 
classical matrix ensembles,  
namely the ones with distribution $\propto e^{-n\tr(V(X))} \ud X$,
for which the deviations are well known (see Theorem \ref{tholdx}). We first 
 tackle a closely related  question:  the 
large deviation properties of the largest eigenvalues of
a deterministic matrix $X_n$  perturbed by a finite rank random matrix.
We show that the law of these extreme eigenvalues satisfies
a large deviation principle for a
 fairly general class of random
finite rank perturbations. We can then consider random matrices $X_n$,
independent of the perturbation, by studying the deviations
of the perturbed matrix conditionally to the non-perturbed matrices.
Even though our rate functions are not very explicit in general, in the simple case
where $X_n=0$, we can retrieve more explicit formulae (see Section
\ref{wishartsec}). In fact, even in this simple case of sample covariance matrices with non-Gaussian entries,
 our large deviation result seems to be new
and improves on \cite{Fey}.

 Our approach is based, as in  \cite{CDF09,CDF09b,bai-yao-TCL},
on the characterization of the eigenvalues 
via the determinant of a matrix with fixed size : it is an $r\times r$ matrix  whose entries 
are the Stieltjes transforms of
the non-deformed matrix evaluated along the random
vectors of the perturbation. We obtain a large deviation principle for the
law of this characteristic polynomial (seen as a continuous function outside
of the spectrum of the deterministic matrix)
by classical large deviation techniques. Even though the application 
which associate to a function its zeroes is not continuous for the
weak topology, we deduce from the latter a large deviation principle for the
law of the zeroes of this characteristic polynomial, that is the extreme eigenvalues of the deformed matrix model.


\section{Statement of the results}

\subsection{The models}
Let $X_n$ be a real diagonal matrix of size $n \times n$ with eigenvalues 
$\la_1^{n} \geq \la_2^{n} \ge \ldots  \ge \la_n^{n}.$

We  perturb $X_n$ by a random matrix whose rank does not depend on $n$.  More precisely, let $m, r$ be fixed positive integers  
and $\theta_1 \ge \theta_2 \ge \ldots \ge \theta_m >0 > \theta_{m+1} \ge \ldots \ge \theta_r$ be fixed,
let   $G= (g_1, \ldots, g_r)$  be a random vector and
$(G(k)=(g_1(k),\ldots,g_r(k)))_{  k\geq 1 }$  be independent copies of $G.$
We then define the $r$
 vectors with dimension $n$
$$G_1^{n}:=(g_1(1),\ldots, g_1(n) )^T\; ,\ldots\ldots \ldots,\; G_r^{n}:=(g_r(1), \ldots,g_r(n) )^T$$  
and study  the eigenvalues $\wtl_1^{n} \ge \cdots \ge \wtl_n^{n}$ of the deformed matrices
\be\label{def_model_iid_29.06}\wtX = X_n + \ff{n}\sum_{i=1}^r \theta_i G_i^{n}{G_i^{n}}^*.\ee
In the sequel, we will refer to  the model  \eqref{def_model_iid_29.06}   as the {\it i.i.d. perturbation model}.\\

Alternatively, if we assume moreover that the law of $G$ does not charge any hyperplane, then,  for $n>r,$
the $r$ vectors $G_1^{n}, \ldots, G_r^{n}$ 
 are almost surely linearly independent
and we denote by $(U_i^{n})_{1 \le i \le r}$ the  vectors
 obtained from  $(G_i^{n})_{1 \le i \le r}$ by a Gram-Schmidt
 orthonormalisation
procedure with respect to the usual scalar product on $\C^n.$
 We shall then consider the eigenvalues $\wtl_1^{n} \ge \cdots \ge \wtl_n^{n}$ of 
\be\label{def_model_ortho_29.06}\wtX = X_n + \sum_{i=1}^r \theta_i U_i^{n}{U_i^{n}}^*\ee 
and  refer in the sequel to the model \eqref{def_model_ortho_29.06} as   the  {\it orthonormalized perturbation model.}

 If $g_1, \ldots, g_r$ are $r$ independent  standard (real or complex) Gaussian variables,  it is well known that the law of $(U_i^{n})_{1\le i\le r}$ 
is the uniform measure on the set of $r$ orthonormal vectors.
The  model  \eqref{def_model_ortho_29.06} coincides then  with the one introduced in \cite{benaych-rao.09}.\\

Our goal will be to examine the large deviations for the $m$ largest eigenvalues of 
the deformed matrix $\wtX,$ with $m$ the number of positive eigenvalues of the random deformation.

\subsection{The assumptions}
Concerning the spectral measure 
of the full rank deterministic matrix $X_n,$
we assume the following

\begin{assum} \label{spacing} 
The empirical distribution  $\frac{1}{n} \sum_{i=1}^n \delta_{\la_i^{n} }$ of $X_n$  converges weakly as $n$ goes to infinity to a
compactly supported \pro $\mu$.
\end{assum}

\noindent
Concerning the random vector $G$, we make  the following assumption. It allows to claim that with \pro one, the column vectors $G_1^n, \ldots, G_r^n$ are linearly independent and   is   technically needed in the proof of Lemma \ref{norm}. It  is also the reason why we say that the column vectors $G_1^n, \ldots, G_r^n$ or $U_1^n, \ldots, U_r^n$ are {\it delocalized} with respect to the eigenvectors of $X_n$. Indeed, the eigenvectors of $X_n$ are the vectors of the canonical basis, whereas we know that with \pro one, none of the  entries of the $G_i^n$'s (or of the $U_i^n$'s) is zero. The i.i.d. feature of the $G(k)$'s allows even to assert that all entries of each $G_i^n$'s (or of the  $U_i^n$'s) have the same distribution. 

\begin{assum} \label{hyponG}
$G=(g_1,\ldots,g_r)$ is a   random vector with entries in $\mathbb{K}=\R$ or $\C$   
such that there exists $\alpha >0$ with
$\E(e^{\alpha  \sum_{i=1}^r |g_i|^2}) < \infty.$ In the
orthonormalized perturbation model, we 
 assume moreover that
 for any $\lambda \in \mathbb{K}^r\bck\{0\}$,  $\Pro(\sum_{i=1}^r \lambda_i g_i =0)=0$ 
\end{assum}

The law of $G$ could also depend on $n$ provided it satisfies
the above hypothesis uniformly on $n$ and converges in law as $n$ goes to infinity. 

We 
 consider two distinct kind of assumptions on the extreme
 eigenvalues of $X_n.$ 
 We will be first interested in the case when these extreme 
eigenvalues  {\bf stick to the bulk} (see Assumption \ref{without}), and then to the case \textbf{with outliers}, when we allow some eigenvalues
of $X_n$ to take their limit outside the support of the limiting measure $\mu$ (see Assumption \ref {with}).

\subsection{The results in the case without outliers}

We first consider the case where the  extreme 
eigenvalues  of $X_n$ {stick to the bulk}.
\begin{assum} \label{without}
 The largest and smallest eigenvalues of $X_n$  tend  
respectively to the upper bound (denoted by 
 $b$) and the lower bound (denoted by 
 $a$) of the support of $\mu$.
\end{assum}

Our main theorem is the following (see 
 Theorem \ref{main} and Theorem \ref{min} for precise statements).

\begin{Th}\label{mainintro}
Under Assumptions \ref{spacing}, \ref{hyponG} and \ref{without}, the law of the $m$ largest eigenvalues $(\wtl_1^{n},\ldots, \wtl_m^{n})
\in\mathbb R^m$
of $\wtX$ satisfies a large deviation principle (LDP) in the scale $n$
with a good rate function $L$. In other words, for any $K\in\mathbb R^+$,
$\{L\le K\}$ is a compact subset of $\mathbb R^m,$
for any closed
set $F$ of $\mathbb R^m$,
$$\limsup_{n\ra\infty}\frac{1}{n}\log \Pro\left((\wtl_1^{n},\ldots, \wtl_m^{n})\in F\right)\le -\inf_F L$$
and for any open set $O\subset\mathbb R^m$,
$$\liminf_{n\ra\infty}\frac{1}{n}\log \Pro\left((\wtl_1^{n},\ldots, \wtl_m^{n})\in O\right)\ge -\inf_O L.$$
 Moreover, this rate function achieves its minimum
value at a unique $m$-tuple $(\lambda_1^*,\ldots,\lambda_m^*)$
towards which $(\wtl_1^{n},\ldots, \wtl_m^{n})$
converges almost surely.
\end{Th}
Theorem  \ref{mainintro} is true for both 
 the i.i.d. perturbation model and the orthonormalized perturbation model, but the exact expression of the rate function $L$ is not the same for both models. As could be expected, the minimum $(\lambda_1^*,\ldots,\lambda_m^*)$
only depends on the $\theta_i$'s, on the limiting spectral distribution $\mu$ of $X_n$, and on the covariance matrix of the vector $G,$
this latter dependence coming from the fact that the rate function involves a Laplace transform of the law of $G$ and its behavior near the extremum will generically be governed by the second derivatives, that is the covariance. \\

The rate function $L$ is not explicit in general. However, in the particular case where $X_n=0$, $L$ can be evaluated. It amounts to consider 
the large deviations of the eigenvalues   of matrices $W_n=\frac{1}{n}G_n^*\Theta G_n$ for $G_n$ an $n\times r$ matrix, with $r$ fixed and $n$ growing to infinity. $L$ is very explicit when $G$ is Gaussian but
even when the entries are not Gaussian, we can
 recover a large deviation principle and 
refine a bound of \cite{Fey}  about the deviations
of the largest eigenvalue (see Section \ref{wishartsec}).

\subsection{The results in the case with outliers}

We now consider the case where some eigenvalues
of $X_n$ escape from the bulk,
so that  Assumption \ref{without} is not fulfilled.
We  assume that these eigenvalues, that we call
{\it outliers}, converge: 

\begin{assum} \mbox{}
 \label{with}
  There exist some non negative integers $p^+, p^- $ such that for any $i \le p^+,$ $\la_i^{n} \ninf\ell_i^+,$ 
for any $j \le p^-,$  $\la_{n-j+1}^{n}\ninf\ell_j^-,$  $\la_{p^++1}^{n} \ninf b$ and   $\la_{n-p^-}^{n} \ninf a$ with 
$-\infty<\ell^-_1 \le \ldots \le \ell_{p^-}^-<a \leq b < \ell_{p^+}^+ \leq \ldots \le \ell_1^+ < \infty,$
where $a$ and $b$ denote respectively the lower and upper bounds of the support of the limiting measure $\mu.$
\end{assum}

To simplify  the notations in the sequel we will use the following conventions : $\ell_{p^-+1}^-:= a$ and $\ell_{p^++1}^+:=b.$

In this framework, we will need to make on $G$ the additional following
assumption.\\

\begin{assum} \mbox{}
 \label{Gwith}
The law of the vector $\frac{ G }{\sqrt{n}}$ 
satisfies a large deviation principle in the scale $n$ with a good rate function that we denote by $I$. 
\end{assum}

\begin{Th}\label{mainout0}
If Assumptions \ref{spacing}, \ref{hyponG}, \ref{with} and \ref{Gwith} hold,
the law of the $m+p^+$ largest eigenvalues of $\wtX$ 
satisfies a large deviation principle with a good rate function $L^o$.
\end{Th}

Again, Theorem  \ref{mainout0} is true for both  i.i.d. perturbation model and orthonormalized perturbation model, but the rate function is not the same for both models. A precise definition of $L^o$ will be given in Theorem \ref{mainout}.

Before going any further, let us discuss Assumption \ref{Gwith}. On one side, let us give some natural examples for which
the assumtion is fulfilled.
\begin{lem}\label{verias}
\begin{enumerate}
\item
If $G=(g_1,\ldots,g_r)$ are i.i.d standard Gaussian
variables, Assumption \ref{Gwith} holds with $I(v)=\frac{1}{2}\|v\|_2^2$.
\item If $G$ is such that for any $\alpha >0,$ $\E[e^{\alpha\sum_{i=1}^r |g_i|^2}]< \infty$, then 
Assumption \ref{Gwith}
 holds with $I$   infinite except at $0$, where it takes value $0$. 
\end{enumerate}
\end{lem}
\begin{pr}
The first result can be seens as  a direct consequence of Schilder's theorem.
For the second, it is enough to notice by Tchebychev's  inequality
that for all $L,\delta>0$,
$$\Pro\left( \max_{1\le i\le r}|g_i|^2\ge \delta n\right)
\le r e^{-L\delta n}\E(e^{L\sum_{i=1}^r  |g_i|^2})$$
so that taking the large $n$ limit and then $L$ going to infinity yields 
for any $\delta>0$
$$\limsup_{n\ra\infty} \frac{1}{n} \log
\Pro\left( \max_{1\le i\le r}|g_i/\sqrt{n}|^2\ge \delta\right)=-\infty$$
thus proving the claim.
\end{pr}\\

On the other side, we want to emphasize that in the case with outliers, the individual LDP stated in Assumption \ref{Gwith}
will be  crucial. To understand more deeply this phenomenon, we refer the interested reader  to some couterexamples
when this assumption is not fulfilled that are studied in \cite[Section 2.3]{jamal-2002} and a related discussion in the introduction
of \cite{mylene-sandrine-jamal}.\\

\subsection{Large deviations for the largest eigenvalues of perturbed matrix models}\mbox{} \\
\label{secmm}

We apply hereafter the results above to study the large deviations
of the law of the extreme eigenvalues of perturbations of randomly chosen matrices $X_n$
distributed according to 
  the Gibbs measure
$$\ud\mu^n_\beta(X)=\frac{1}{Z_n^\beta}e^{-n\Tr(V(X))} \ud^\beta X$$
with $\ud^\beta X$ the Lebesgue measure on the set of $n\times n$ Hermitian matrices if $\beta=2$ (corresponding to  $G$  $\C^r$-valued)
or $n\times n$ symmetric matrices  if $\beta=1$ (corresponding to $G$ $\R^r$-valued). \\

Let us first recall a few facts about the non-perturbed model.
It is well known that if $X_n$ is distributed according to $\mu^n_\beta,$ the law of the eigenvalues of $X_n$
is given by 
$$\Pro_{V,\beta}^n(\ud\lambda_1,\ldots,\ud\lambda_n)=\frac{\mathds 1_{\lambda_1>\lambda_2>\cdots>\lambda_n}
}{Z^n_{V,\beta}}\prod_{1\le i<j\le n}|\lambda_i-\lambda_j|^\beta
e^{-n\sum_{i=1}^n V(\lambda_i)}\prod_{i=1}^n \ud\lambda_i.$$

We will make on the potential $V$ the following assumptions :
\begin{assum}\label{assX}
\begin{itemize}
\item[$i)$] $V$ is continuous with values in $\R\cup\{+\infty\}$ and  
$$\liminf_{|x|\ra\infty} \frac{V(x)}{\beta\log |x|} >1.$$
\item[$ii)$] For all integer numbers $p$, the limit
$$\lim_{n\ra\infty}\frac{1}{n}\log \frac{Z^{n-p}_{nV/n-p,\beta}}{Z^n_{V,\beta}}
$$
exists and is denoted by $\alpha_{V,\beta}^p$.
\item[$iii)$] Under $\Pro_{V,\beta}^n,$ the largest eigenvalue $\lambda_1^{n}$
converges  almost surely 
to the upper  boundary $b_V$ of the support of  $\mu_V$.
\end{itemize}
\end{assum} 

Under part $i)$ of the assumption, one can get a large deviation
principle in the scale $n^2$ for the law 
of the  spectral measure $n^{-1}\sum_{i=1}^n \delta_{\lambda_i}$
under $ \Pro_{V,\beta}^n$
(see \cite{BAG3}), resulting in particular with
the almost sure convergence of the spectral measure to a probability 
measure $\mu_V^\beta$. If we add part $ii)$ and $iii),$ one can derive 
the large deviations for the extreme eigenvalues of $X_n$ 
(see  \cite{bdg01}, and also \cite[Section 2.6.2]{alice-greg-ofer}\footnote{Note that in the published version
of \cite{alice-greg-ofer},
part $iii)$ was not mentioned but it appears in the errata sheet available online : \url{http://www.wisdom.weizmann.ac.il/~zeitouni/cormat.pdf}}). We give below a slightly more general statement to
consider the deviations of the $p$th largest eigenvalues (note that
the $p$th smallest can be considered similarly). 


One can notice that these assumptions hold in a wide generality.
In particular, they are satisfied for the law of the GUE ($\beta=2$, $V(x)=x^2$)
and the GOE ($\beta=1$, $V(x)=x^2/2$) as part $ii)$ is verified by Selberg formula 
whereas part $iii)$ is well known (see \cite[Section 2.1.6]{alice-greg-ofer}).
For the case of Gaussian Wishart matrices, we know (see e.g. \cite[p 190]{alice-greg-ofer}) that the joint law of the eigenvalues can be written 
as $ \Pro_{V_{p,n},\beta}^n$ with $V_{p,n}(x) = \frac{\beta}{4} x - (\beta[1-\frac{p}{n}+\frac{1}{n}]-\frac{1}{n})\log x$ on $(0,\infty).$
If the ratio $\frac{p}{n}$ converges to $\alpha,$ one can easily show that the law of the  largest eigenvalues
are exponentially equivalent under $ \Pro_{V_{p,n},\beta}^n$  and under $ \Pro_{V,\beta}^n,$ with $V(x) =  \frac{\beta}{4} x - \beta(1-\alpha)\log x$   on $(0,\infty),$ for which the assumptions are satisfied.\\

\begin{Th}\label{tholdx}
Under Assumption \ref{assX}, 
the law of the  $p$ largest eigenvalues $(\lambda_1^{n}>\cdots>
\lambda_p^{n})$  of $X_n$
satisfies a large deviation principle in the scale $n$
and with good rate function given by
$$J^p(x_1,\ldots,x_p)= \left\{ \begin{array}{ll}
\sum_{i=1}^p J_V(x_i)+p \alpha_{V,\beta}^1, & \textrm{if } x_1\ge x_2\ge \cdots \ge x_p, \\ \\
\infty, & \textrm{otherwise,}\end{array}
                                                                                         \right.$$
with $J_V(x)=V(x)-\beta\int \log|x-y|d\mu_V(y)$.
\end{Th}
\begin{rmq}
Note that in the case of the GOE and the GUE (see \cite{bdg01}),
$$J_V(x)=\beta\int_2^x\sqrt{(y/2)^2-1} dy-\alpha^1_{V,\beta},\quad \alpha^1_{V,\beta}=-{\beta}/2.$$
\end{rmq}

Let us now go to the perturbed model. An important remark is that,
due to the rotational invariance of the law of $X_n,$ one can in fact consider very general orthonormal perturbations.
We make the following 

\begin{assum}\label{mmperturb}
 $(U^n_1, \ldots, U^n_r)$ is a family of orthonormal vectors in $(\mathbb R^n)^r$ (resp. $(\mathbb C^n)^r$) if $\beta= 1$
(resp. $\beta= 2$),
either deterministic
or independent of $X_n.$ 
\end{assum}

Indeed, 
under these assumptions,  $\wtX$ has in law the same eigenvalues as \linebreak
$ D_n + \sum_{i=1}^r \theta_i (O_nU^n_i)(O_nU^n_i)^*,$ with $D_n$ a real diagonal matrix with $\Pro_{V,\beta}^n$-distributed eigenvalues   and $O_n$  Haar distributed on the orthogonal (resp. unitary) group of size $n$ if $\beta= 1$
(resp. $\beta= 2$), independent of $\{D_n\}\cup\{U^n_1, \ldots, U_r^n\}$. Now, from the well know properties of the Haar measure, if the $U^n_i$'s satisfy Assumption \ref{mmperturb},
then the $O_nU^n_i$'s are column vectors of a  Haar distributed matrix. In particular  they can be obtained by the orthonormalization procedure described in the introduction, with $G = (g_1, \ldots,g_r)$ a vector whose components are  i.i.d. Gaussian standard variables (which satisfies in particular
Assumption \ref{Gwith}).

With these considerations in mind, we can state the large 
deviation principle for 
the extreme eigenvalues of $\wtX$. We recall that $b_V$ is the rightmost point of the support of $\mu_V$. 

\begin{Th}\label{theorandom}
With $V$ satisfying Assumption \ref{assX}, we consider the orthonormalized perturbation model
under Assumption  \ref{mmperturb}. Then, for any integer $k,$ 
the law of the  $k$ largest eigenvalues $(\wtl_1^{n},\cdots,
\wtl_k^{n})$  of $\wtX$
satisfies a large deviation principle in the scale $n$
and with good rate function given by
$$\tilde J^k(x_1,\ldots,x_k)= 
\inf_{p\ge 0}\inf_{\ell_1\ge \cdots\ge\ell_p>
b_V} \{ L^0_{\ell_1,\ldots,\ell_p}
(x_1,\ldots, x_k) + J^p(\ell_1,\ldots,\ell_p)\},$$
if $x_1\ge \cdots\ge x_k$, the function being infinite otherwise.\\                                                                             
Here, $L^0_{\ell_1,\ldots,\ell_p}$ is the rate function defined in
Theorem \ref{mainout} for the orthonormalized perturbation model built on $G=(g_1, \ldots,g_r)$ i.i.d. standard Gaussian variables
and $X_n$ with limiting spectral measure $\mu_{V}$ and outliers  $\ell_1,\ldots,\ell_p.$
\end{Th}

\section{Scheme of the proofs}\label{Sectionscheme}
 
The strategy of the proof will be quite similar  in both  cases (with or without outliers), so,  for the
 sake of simplicity, we will 
outline it in the present section only in the  case without outliers  (both 
 the i.i.d. perturbation model and the orthonormalized perturbation model will be treated simultaneously).

The cornerstone  is 
 a nice representation, already crucially used  in many
papers on finite rank deformations (see e.g.
 \cite{benaych-rao.09,bai}), of the eigenvalues 
$(\wtl_1^{n},\ldots, \wtl_m^{n})$
as zeroes of a fixed deterministic polynomial in the entries
 of matrices of  size $r$ 
depending only on the resolvent of $X_n$
and the random vectors $(G^n_i)_{1\le i\le r}$.

Indeed, if $V$ is the $n\times r$  matrix with column vectors $\bbm U_1^{n}\cdots 
U_r^{n}\ebm$ in the orthonormalized perturbation model and  $\bbm G_1^{n}\cdots 
G_r^{n}\ebm$ in the i.i.d. perturbation model, $\Theta$ the matrix $\diag(\theta_1, \ldots, \theta_r)$ and $I_n$ the identity in $n\times n$ matrices, 
 the characteristic polynomial 
 of $\wtX$  reads 
\begin{equation}\label{51}
\det(zI_n-\wtX)=\det(zI_n-X_n-V \Theta V^*)
=\det(zI_n-X_n)\det(I_r-V^*(zI_n-X_n)^{-1}V\Theta)
\end{equation}

It means that the eigenvalues of $\wtX$ that are not\footnote{ We show 
in section \ref{SectionAppendixdistinct} that  the spectra of $X_n$
and $\wtX$ are disjoint in generic situation.} 
 eigenvalues of $X_n$ 
are the zeroes of  $\det(I_r-V^*(zI_n-X_n)^{-1}V\Theta),$ which is the determinant of
a matrix whose size is independent of  $n$. 

Because of the relation between  $V$ and  the random vectors 
$G_1^{n},\ldots, G^{n}_r$, it is not hard to check that, if we let, for  
$z \notin \{\la_1^{n}, \ldots, \la_n^{n}\},$  $K^{n}(z)$  and $C^{n}$ be the elements of the set $\Hr$  of $r\times r$ Hermitian matrices given, for $1\le i \le j \le r$, by
\be \label{defK} K^{n}(z)_{ij} = \frac{1}{n}
 \sum_{k=1}^n \frac{\overline{g_i(k)}g_j(k)}{z-\la_k^{n}}\,\ee
and \be \label{defC} C^{n}_{ij}=\frac{1}{n}\sum_{k=1}^n \overline {g_i(k)}g_j(k)\,,\ee
we have (see Section \ref{secpol} for details):
\begin{propo}\label{zeroesscheme}
In both i.i.d and orthonormalized perturbation models,
 there exists a function
 $P_{\Theta, r} $ defined on  ${\Hr}\times {\Hr}$
which is polynomial in the entries of its arguments and 
 depends only on the matrix $\Theta,$ such that  any $z \notin \{\la_1^{n}, \ldots, \la_n^{n}\}$
is an eigenvalue of $\wtX$ if and only if 
$$ H^{n}(z):=P_{\Theta, r}(K^{n}(z), C^{n})=0\,.$$

\end{propo} 
Of course, the polynomial $P_{\Theta, r}$
is different in the i.i.d. perturbation model and the orthonormalized  perturbation model. In the i.i.d. perturbation model,   $P_{\Theta, r} $ is simpler 
and does not depend on $C$.
This proposition characterizes the eigenvalues of $\wtX$ 
as the zeroes of the random function $H^{n}$, which depends continuously
(as a polynomial function) on the random pair $ (K^{n}(\cdot), C^{n})$.
The large deviations of these eigenvalues are therefore
inherited  from the large deviations of $ (K^{n}(\cdot), C^{n})$, which we thus study in detail before getting into the deviations of the eigenvalues themselves.
Because  $K^{n}(z)$
blows up when $z$ approaches $\lambda_1^{n}$, which itself converges 
to $b$, we  study the large deviations of $ (K^{n}(z), C^{n})$ for $z$ away from $b$.
We shall let  $\K$ be a  compact interval 
 in $(b, \infty),$  $\mathcal C(\K, {\Hr})$ and  $\mathcal C(\K,\mathbb R )$ be  the space of continuous functions on $\K$
taking values respectively in ${\Hr}$ and in $\mathbb R.$ We endow
the latter set with the uniform
topology.  We will then prove that (see Theorem \ref{theo1} for a precise statement and 
a definition of the rate function {\bf I} involved)

\begin{propo}\label{theo1scheme}
The law of $ ((K^{n}(z))_{z\in \K}, C^{n})$ on $ \mathcal C(\K, {\Hr})\times {\Hr}$ equipped with the uniform topology,
satisfies a large deviation principle in the scale $n$ and with
good rate function ${\bf I}$.
\end{propo}
By the contraction principle,
we therefore deduce 
\begin{cor}
The law of $(H^{n}(z))_{z \in \K}$ on  $\mathcal C(\K,\mathbb R )$  equipped with the uniform topology,
satisfies a large deviation
principle in the scale $n$ and with rate function
given, for a continuous function $f\in \mathcal C(\K,\mathbb R)$,
by 
$$J_\K(f)=\inf\{ \mathbf I(K(\cdot), C)\ste(K(\cdot), C)\in \mathcal C(\K, {\Hr})\times {\Hr},
P_{\Theta,r}(K(z), C))=f(z)\,\,
\forall z \in \K\}
$$
with $P_{\Theta,r}$ the polynomial function of Proposition \ref{zeroesscheme}.
\end{cor}
Theorem \ref{mainintro}
is then a consequence of this corollary with,
heuristically,  $L(\alpha)$ the infimum of $J_{[b,+\infty)}$
on the set of functions which vanish 
exactly  at $\alpha \in \R^m$. An important technical issue will come from  the
fact that the set of functions which vanish exactly  at $\alpha$
has an empty interior, which requires extra care for
the large deviation lower bound.\\

The organisation of the paper will follow the scheme we have just  described: in the next section, we detail the 
orthonormalization procedure and prove Proposition \ref{zeroesscheme}. 
Section \ref{SectionLDPH} and Section \ref{SectionLDPvap} will then deal more specifically with the  case without outliers.
In Section \ref{SectionLDPH}, we  establish
the functional large deviation principles for $ (K^{n}(\cdot), C^{n})$ and $H^{n}$, whereas  
  Section \ref{SectionLDPvap} is devoted to the proof of 
 our main results in this case,
namely the large deviation principle for the largest eigenvalues of $\wtX$ and the almost sure convergence to the minimisers
of the rate function. In Section \ref{wishartsec}, we will see that the rate function can be studied further in the special case when 
$X_n = 0.$  We then turn to the case with outliers in Sections \ref{section_main_outliers_1_30.6.10} and \ref{section_main_outliers_30.6.10}. Therein, the proofs will be less detailed,
but we will  insist on the points that differ from the previous case. 
The extension to random matrices
$X_n$  given by  classical matrix models is presented  in Section \ref{Proof_classical_300610}.
To make the core of the paper easier to read, we gather some technical results in Section \ref{SectionAppendix}.\\

\section{Characterisation of the eigenvalues of $\wtX$ as   zeroes of a function $H^n$}\label{secpol}

The goal of this section is to prove Proposition \ref{zeroesscheme}.
As will be seen further, the proof of this proposition is straightforward in the i.i.d. perturbation model but more involved in the orthonormalized
perturbation model and we first detail the orthonormalization procedure.
 
\subsection{The Gram-Schmidt orthonormalisation procedure} \label{SectionGS}\mbox{}\\

We start  by detailing the construction of $(U_i^{n})_{1\le i\le r}$ from $(G_i^{n})_{1\le i\le r}$ in the orthonormalized perturbation model. 
The canonical scalar product in $\C^n$ will be denoted by $\langle v, w\rangle =v^*w= \sum_{k=1}^n \overline{v_k}w_k$,
 and the associated norm  by $\|\cdot \|_2.$  We also recall that ${\Hr}$ is the space of $r\times r$ 
either symmetric or Hermitian matrices, according to whether $G$ is a real ($\mathbb{K}=\R$) or complex ($\mathbb{K}=\C$) random vector.

Fix $1\le r\le n$ and  consider a linearly independent family $G_1, \ldots, G_r$ of vectors in $\mathbb{K}^n$. 
Define their  Gram matrix (up to a factor $n$) 
$$C=[C_{ij}]_{i,j=1}^r,\quad \textrm{ with } C_{ij}=\ff{n}\lan G_i,G_j\ran.$$ 
We then define
  \begin{equation}\label{defq}
 q_1= 1 \qquad\textrm{ and \qquad for } i=2, \ldots, r, \quad q_i:=\det [C_{kl}]_{k,l=1}^{i-1}\end{equation}
and the  lower triangular matrix  $A=[A_{ij}]_{1\leq j\le i\le r}$  as follows~:
for all $1\leq j < i\leq r$, 
\be \label{Aij}
A_{ij}=\f{\det [\gamma_{k,l}^j]_{k,l=1}^{i-1}}{q_i} \textrm{ with } 
\gamma_{kl}^j = \left\{ \begin{array}{ll}
                   C_{kl},& \textrm{si } l\neq j \\
- C_{ki},& \textrm{si } l= j .
                  \end{array}
 \right.
\ee
Note that by linear independence of the $G_i$'s, none of the $q_i$'s is zero
so that the matrix $A$ is well defined.

Then the vectors $W_1, \ldots, W_r$ defined, for $i=1, \ldots, r$, by $$W_i=\sum_{l=1}^iA_{il}\f{G_l}{\sqrt{n}} $$ are  orthogonal and the $U_i$'s, defined, for $i=1, \ldots, r$, by $$U_i=\f{W_i}{\|W_i\|_2}$$ are  
orthonormal. They are said to be the Gram-Schmidt orthonormalized vectors from $(G_1, \ldots, G_r).$
The following proposition, which can be easily deduced from the definitions we have just introduced, will be useful in the sequel.  

\begin{propr}\label{pol}
For each $i_0=1, \ldots, r$, there is a real function $P_{i_0}$, defined on ${\Hr}$,
polynomial in the entries of the matrix, not depending on $n$ and nor on the $G_i$'s,  \st $$\|q_{i_0}W_{i_0}\|^2_2=P_{i_0}(C).$$ 
Moreover, the polynomial function $P_{i_0}$ is positive on the set of positive 
definite matrices.
\end{propr}

The last assertion of the proposition comes from the fact that any positive definite $r\times r$ Hermitian matrix  is the Gram matrix of a linearly independent family
 of $r$ vectors of $\mathbb{K}^r$ (namely the columns of its square root).\\

Let now $G$ be a random vector satisfying  Assumption \ref{hyponG} and  $(G(k), k\ge 1)$ be i.i.d. copies of $G$.  Let $G_i^{n}=(G(k)_i)_{1\le k\le n}$ for $i\in\{1,\ldots,r\}$. One can easily check that if $n >r$, these vectors   are almost surely 
linearly  independent, so that we can 
apply Gram-Schmidt orthonormalisation
 to this family  of random vectors. We define the $r\times r$ matrices $C^{n}, A^{n}$, the real number $q_i^{n}$ 
and the vectors $W_1^{n}, \ldots, W_r^{n}$, $U_1^{n}, \ldots, U_r^{n}$  of $\mathbb{K}^n$  as above.
As announced in Section \ref{SectionDefmodel}, these $U_i^{n}$'s are the Gram-Schmidt orthonormalized  
of the $G_i^{n}$'s we used to define our model in the introduction.\\

\subsection{Characterization of the eigenvalues of  $\wtX$: proof of Proposition \ref{zeroesscheme}}\mbox{}\\

As explained in Section  \ref{Sectionscheme}, a crucial observation (see \cite[Proposition 5.1]{benaych-rao.09}) is that the eigenvalues of  $\wtX$ can be characterized as the zeroes of 
a polynomial function of matrices of size $r \times r.$  This was stated in  Proposition \ref{zeroesscheme} which we prove below. 

\noindent{\it Proof of Proposition \ref{zeroesscheme}.} We first recall \eqref{51}, that is for $z \notin \{\la_1^{n}, \ldots, \la_n^{n}\},$
\beq
\det(zI_n-\wtX)&=&\det(zI_n-X_n-V \Theta V^*)\\
&=&\det(zI_n-X_n)\det(\Theta)\det(\Theta^{-1}-V^*(zI_n-X_n)^{-1}V)
\eeq
Hence any  $z \notin \{\la_1^{n}, \ldots, \la_n^{n}\}$
is an eigenvalue of $\wtX$ if and only if 
$$ D_n(z) := \det(\Theta^{-1}-V^*(zI_n-X_n)^{-1}V) = 0.$$

We denote  by $\mathbf G $  the $n \times r$ matrix with column vectors 
$(G_i^{n})_{1 \le i \le r},$ so that   $K^{n}(z) = \frac{1}{n} \mathbf G^* (zI_n-X_n)^{-1}\mathbf G.$

In the i.i.d. perturbation model, as $V= \mathbf G,$ Proposition \ref{zeroesscheme} follows immediately with
 $$H^n(z):=\det(\Theta^{-1}-V^*(zI_n-X_n)^{-1}V),$$ which is actually a polynomial, depending on $\Theta$, in the entries of $K^n(z)$.

In the orthonormalized perturbation model, the Gram-Schmidt procedure makes things a bit more involved. 

If we denote by $D$ the $r \times r$ diagonal matrix given by $D = \diag (\|W_1^{n}\|_2, \ldots, \|W_r^{n}\|_2) $
and $\Sigma = (A^{n})^T,$ then $V $ is equal to
$ n^{-1/2}\mathbf G \Sigma D^{-1}$
and  we deduce that
$$  D_n(z) = \det (\Theta^{-1} -  D^{-1}\Sigma^*K_n(z)\Sigma D^{-1} ).$$

Now, if we define  $Q = \diag (q_1^{n}, \ldots, q_r^{n})$ (recall \eqref{defq}), $E= DQ$, $F=\Sigma Q$ and
$H^{n}(z) := \det(E^*\Theta^{-1}E -F^* K_n(z)F)$ then on one hand, one can check that
 $$D_n(z) = (\det E^*E)^{-1}   H^{n}(z),$$ so that 
any  $z \notin \{\la_1^{n}, \ldots, \la_n^{n}\}$
is an eigenvalue of $\wtX$ if and only if it is a zero of $H^{n}.$\\
On the other hand, 
$H^{n}(z)$ is obviously a polynomial (depending only on the matrix $\Theta$) of the entries of $K^{n}(z)$, $E^*\Theta^{-1}E$ and $F$. 
Furthermore, $E^*\Theta^{-1}E$ is a diagonal matrix whose $i$-th entry is given by 
$(E^*\Theta^{-1}E)_i = \theta_i^{-1} \|q_i^{n} W_i^{n}\|^2_2= \theta_i^{-1} P_i(C^{n})$ (by Property \ref{pol}) 
and   $F_{ij} =\det [\gamma_{k,l}^j]_{k,l=1}^{i-1}$ with $\gamma_{k,l}^j$ defined in \eqref{Aij}. This concludes the proof.\hfill$\square$

\section{Large deviations   for   $H^{n}$ in the  case without outliers}
\label{SectionLDPH}
 \noindent
We assume throughout this section that Assumptions \ref{spacing}, \ref{hyponG} and \ref{without} hold.

\subsection{Statement of the result}\mbox{}\\

 In the sequel, $\K$ will denote any compact interval
included in $(b, \infty),$ and we denote by $z^*$ its upper bound. 
We equip $\mathcal C(\K, {\Hr})\times {\Hr}$ with the uniform topology which is given by the distance $d$ defined, for $(K_1,C_1), (K_2,C_2) \in \mathcal C(\K, {\Hr})\times {\Hr}$ by
$$ d((K_1,C_1), (K_2,C_2)) = \sup_{z \in \K} \|K_1(z) - K_2(z)\|_2 +  \|C_1 - C_2\|_2,$$ where $\|M \|_2=\sqrt{\Tr(M^2)}$ for all   $M\in {\Hr}.$ 

With $G=(g_1, \ldots,g_r)$ satisfying Assumption \ref{hyponG}, we define $Z$ a matrix in ${\Hr}$ such that, for $i \leq j,$ $Z_{ij}= \overline{g_i}g_j$
and  $\Lambda$ given, for any $H \in {\Hr}$ by
\be\label{1.02.10.def_La} \Lambda(H) = \log \E\left(e^{\Tr(HZ)}\right).\ee

The goal of this section is to show the following theorem.
\begin{Th}\label{theo1}
\begin{enumerate}
\item  \label{theoK} The law of 
 $\left((K^{n}(z))_{z \in \K}, C^{n}\right)$,
viewed as an element of
the space $\mathcal C(\K, {\Hr})\times {\Hr}$ 
equipped with the uniform topology,
satisfies a 
large deviation
principle in the scale $n$ and with good rate function $\bf I$
which is infinite if $K$ is not Lipschitz continuous and otherwise
defined, for $K \in \mathcal C(\K, {\Hr})$ and $C \in  {\Hr}$,
by
 $$\mathbf I(K(\cdot),C)=\sup_{P,X,Y}\left\{ \Tr\left(\int K^\prime(z)
P(z)\ud z+ K(z^*) X+ C Y\right)
-\tilde\Gamma(P,Y,X)\right\}
$$
where  $\tilde\Gamma(P, Y,X)$ is given by the formula
$$\tilde\Gamma(P, Y,X)=\int\Lambda\left(-\int \frac{1}{(z-x)^2}P(z)\ud z+ \frac{1}{z^*-x}X+Y\right)\ud\mu(x)$$
and the supremum is taken 
over
 piecewise constant functions
 $P $  with values in ${\Hr}$ and  $X,Y$ in  ${\Hr}.$

\item  \label{theoH}
The law of $(H^{n}(z))_{z \in \K}$ on   $\mathcal C(\K,\mathbb R)$ 
equipped with the uniform topology,
satisfies a large deviation
principle in the scale $n$ and with rate function
given, for a continuous function $f\in \mathcal C(\K,\mathbb R)$,
by 
$$J_\K(f)=\inf\{ \mathbf I(K(\cdot), C)\ste(K(\cdot), C)\in \mathcal C(\K, {\Hr})\times {\Hr},
P_{\Theta,r}(K(z), C))=f(z)\,\,
\forall z \in \K\}
$$with $P_{\Theta,r}$ the polynomial function of Proposition \ref{zeroesscheme}.

\end{enumerate}
\end{Th}

Since 
 the map $(K(\cdot), C)\longmapsto (P_{\Theta,r}(K(z), C))_{z\in \K}$ from $\mathcal C(\K, {\Hr})\times {\Hr}$ to  $\mathcal C(\K,\mathbb R)$, both equipped 
with their uniform topology, is continuous and $\mathbf I$ is a good rate function,
the second part of the theorem is a direct consequence of its first part and the contraction principle \cite[Theorem 4.2.1]{DZ}.

The reminder of the section will be devoted to the proof of the first part of the theorem and the study of the properties of the rate function
$\bf I,$ in particular its minimisers.

\subsection{Proof of Theorem \ref{theo1}.}\mbox{}\\
\label{subsec:proofhn}

The strategy will be to 
establish a LDP for finite dimensional marginals of
the process $\left((K^{n}(z))_{z\in \K}, C^{n}\right)$
based on \cite[Theorem 2.2]{jamal-2002}  (see also \cite{bdg01} and \cite{BGR97}). From that, we will establish a LDP in the topology of pointwise convergence 
via  the Dawson-G\"artner theorem. As $\left((K^{n}(z))_{z\in \K}, C^{n}\right)$ will be shown to be exponentially tight 
for the uniform topology, the LDP will also hold in this latter topology.

\subsubsection{Exponential tightness}
\label{subsubexp}
We start with the exponential tightness, stated in the following lemma.
As $\K$ is a compact subset of $(b, \infty)$ and the largest eigenvalue $\la_1^{n}$ tends to $b,$ there exists $0<\e <1$
(depending only on $\K$)
such that for $n$ large enough, for any $z \in \K$ and $1\leq i\leq n,$ $z- \la_i^{n} > \e.$ We fix hereafter such an $\e.$ 

 For any $L >0,$ we  define
$$ \mc{C}_{\K,L} :=   \left\{(K,C)\in \mc{C}(\K, {\Hr})\times {\Hr}\ste  \sup_{z \in \K} \|K(z)\|_2 +  \|C\|_2 \le L, \,\, \textrm{$K$ is $\f{L}{2\eps}$-Lipschitz}
\right\}.$$
We have 

\begin{lem} \label{lemtight}
$$ \limsup_{L \ra \infty} \limsup_{n \ra \infty} \frac 1 n \log \Pro\left(((K^{n}(z))_{z\in \K}, C^{n}) \in \mc{C}_{\K,L}^c\right)= -\infty.$$ 
In particular, the law of $\left((K^{n}(z))_{z\in \K}, C^{n}\right)$ is exponentially tight for the uniform topology
on $\mc{C}(\K, {\Hr})\times {\Hr}.$ 

\end{lem}

\begin{pr}
We claim that 
$$ \left\{ \max_{1 \le i \le r} C^{n}_{ii} \leq \frac{\e L}{2r}\right\}
\subset  \left\{ \left((K^{n}(z))_{z\in \K}, C^{n}\right) \in  \mc{C}_{\K,L}\right\}.$$
 Indeed, for $n $ large enough,
$$
|K^{n}(z)_{ij}-K^{n}(z')_{ij}|\le\sqrt{C^{n}_{ii}C^{n}_{jj}}
\frac{|z-z'|}{\e^2},$$
  whereas since $|C^n_{ij}|^2\le C^n_{ii}C^n_{jj}$,
$\| C^{n}\|_2 \le  r\max_{1 \le i \le r} C^{n}_{ii} $
and  $\|K^{n}(z) \|_2 \le \frac{1}{\e} r \max_{1 \le i \le r} C^{n}_{ii} .$\\

Now, by Assumption \ref{hyponG}, let $\alpha >0$ be such that $C:= \E\left( 
e^{\alpha \sum_{i=1}^r  |g_i|^2}\right)<\infty.$
\beqy
 \Pro\left(  \max_{1 \le i \le r} C^{n}_{ii}  >  \frac{\e L}{2r}\right) &\le& r  \Pro\left(   C^{n}_{11}  >  \frac{\e L}{2r}\right) \\
 &\le& r \E\left( e^{\alpha \sum_k |G_1^{n}(k)|^2}\right)e^{- n\alpha  \frac{\e L}{2r}} 
\leq  r  C^n e^{- n\alpha  \frac{\e L}{2r}} \le  e^{- n\alpha  \frac{\e L}{4r}},
\eeqy
where the last inequality holds for $n$ and $L$ large enough. 
This gives
$$ \limsup_{L \ra \infty} \limsup_{n \ra \infty} \frac 1 n \log \Pro\left(((K^{n}(z))_{z\in \K}, C^{n}) \in  \mc{C}_{\K,L}^c\right)= -\infty.$$ 
By the  Arzela-Ascoli theorem,  $ \mc{C}_{\K,L}$ is a compact subset of  $\mathcal C(\K, {\Hr})\times {\Hr}$ for any $L>0$,  
from which we get immediately the second part of the lemma. \end{pr}

\subsubsection{Large deviation principle for finite dimensional marginals}
\label{subsubmargin}
We now  study  the finite dimensional marginals of our process. 
 More precisely, we intend to show the following:
\begin{propo}\label{PGDmarginals}
 Let $M$ be a positive integer and $b< z_1 < z_2 < \cdots < z_M.$ \\
The law of $\left((K^{n}(z_i))_{1\le i \le M}, C^{n}\right)$ viewed as an element of $\Hr^{M+1}$ satisfies a
large deviation principle in the scale $n$ with good rate function $I_M^{z_1, \ldots, z_M}$ defined, 
for $K_1, \ldots, K_M, C \in {\Hr}$ by
$$  I_M^{z_1, \ldots, z_M}(K_1, \ldots, K_M, C) = \sup_{\Xi_1, \ldots, \Xi_M,Y \in {\Hr}} 
\left\{\Tr\left(\sum_{l=1}^M \Xi_l K_l + YC\right)  - \Gamma_M(\Xi_1, \ldots, \Xi_M,Y)\right\},$$ 
with $\Gamma_M(\Xi_1, \ldots, \Xi_M,Y)$ defined by the formula
$$\Gamma_M(\Xi_1, \ldots, \Xi_M,Y) = \int \Lambda \left(\sum_{l=1}^M \frac{1}{z_l-x}\Xi_l +Y\right) d\mu(x),$$
   $\La$ being  given by  \eqref{1.02.10.def_La}.
\end{propo}

\begin{pr}
The proof of the proposition is a direct consequence of Theorem 2.2 of \cite{jamal-2002}. 
Indeed,  let $Z_1$ be the ${\Hr}$-valued random variable
such that for all $1 \le i,j \le r,$
$$ (Z_1)_{ij} = \overline{g_i(1)}g_j(1)$$
and we define $f$ the matrix-valued continuous function with values in $\R^{[(M+1)r]\times r}$ such that, if we denote
by $I_r$ the identity matrix in ${\Hr},$
$$ f(x) = \left( \begin{array}{c}
                  \frac{1}{z_1-x}I_r \\
\vdots \\
\frac{1}{z_M-x}I_r \\
I_r
                 \end{array}
\right).$$
Now, if $(Z_k)_{1 \le k \le n}$ are iid copies of $Z_1,$  we denote by 
$$ L_n := \frac{1}{n} \sum_{k=1}^n f(\la_k^{n})\cdot Z_k = \left( \begin{array}{c}
                 K^{n}(z_1) \\
\vdots \\
K^{n}(z_M) \\
C^{n}
                 \end{array}
\right).$$
A slight problem is that $\frac 1 n \sum_{i=1}^n \delta_{\lambda_i^n}$ do not fulfill Assumption A.1 in \cite{jamal-2002} in the sense 
that this assumption requires that for all $i,$ $\lambda_i^n$ belongs to the support of the limiting measure $\mu.$ Nevertheless,
it is easy to construct (as was done in the proof of Theorem 3.2 in \cite{mylene-sandrine-jamal}) a sequence $\bar \lambda_i^n$
such that $\frac 1 n \sum_{i=1}^n \delta_{\bar \lambda_i^n}$ fulfills Assumption A.1 in \cite{jamal-2002} and  $\bar L_n := \frac{1}{n} \sum_{k=1}^n f(\bar \la_k^{n})$ is exponentially equivalent to $L_n.$
Then from Theorem 2.2 of \cite{jamal-2002}, we get that $L_n$ satisfies an LDP in the scale $n$ with good rate function  $I_M^{z_1, \ldots, z_M}.$ \end{pr}

\subsubsection{Large deviation principle for the law of $((K_n(z))_{z\in {\mathcal K}}, C^{n})$}
\label{subsubldpK}

The next step is to establish a LDP for the law of $\left((K^{n}(z))_{z\in \K }, C^{n}\right)$
associated with  the topology of pointwise convergence.
The following proposition will be a straightforward application of the Dawson-G\"artner theorem on projective limits.

\begin{propo} \label{DG}
  The law of $\left((K^{n}(z))_{z\in \K}, C^{n}\right)$
as an element of  $\mathcal C(\K, {\Hr})\times {\Hr}$ equipped with the topology of pointwise convergence
satisfies a LDP in the scale $n$ with good rate function $\bf J$ defined as follows : for $K \in \mathcal C(\K, {\Hr})$
and $C \in {\Hr},$
$$ \mathbf J(K, C) = \sup_M \sup_{z_1<\cdots< z_M, z_i \in \K} I_M^{z_1,\ldots,z_M} (K(z_1), \ldots ,K(z_M),C).$$
Moreover $ \mathbf J$ equals the rate function $\mathbf I$ given
in Theorem \ref{theo1}.(1).
\end{propo}

\begin{pr}
Let $\mathcal J$ be the collection of all finite subsets of $\K$ ordered by inclusion.
For $j = \{z_1,\ldots,z_{|j|}\} \in \mathcal J$ and $f$ a measurable function from $\K$ to ${\Hr},$
$p_j(f) = (f(z_1), \ldots, f(z_{|j|})) \in \Hr^{|j|}$.\\
We know from Proposition \ref{PGDmarginals} that the law of $(p_j(K^{n}),C^{n})$ satisfies a LDP
with good rate function $I_{|j|}^{z_1, \ldots, z_{|j|}}.$ Moreover, one can check that
the projective limit of the family $\Hr^{|j|} \times {\Hr}$ is $\Hr^{\K} \times {\Hr}$
equipped with the topology of pointwise convergence.\\
Therefore, the Dawson-G\"artner theorem \cite[Theorem 4.6.1]{DZ} proves
the LDP with rate function ${\mathbf J}$.
The identification of ${\mathbf J}$ as $\mathbf I$
is straightforward
as by  a simple change of variables, $ \mathbf J$ is the supremum
of
$$J(\Xi,M,z):=\Tr\left( \sum_{l=1}^{M-1} \Xi(z_l)(K(z_{l+1})-K(z_{l})) +K(z_M)\Xi(z_M)+CY\right)\qquad$$
$$\qquad\qquad\qquad
-\int\Lambda\left(
\sum_{l=1}^{M-1} \Xi(z_l)\left(\frac{1}{z_{l+1}-x}-\frac{1}{z_l-x}\right)+ \frac{\Xi(z_M)}{z_M-x}+Y\right)d\mu(x)$$
over the choices of $\Xi,M,z$. We may assume
without loss of generality that $z_M=z^*$.
 Putting $P(z)=\sum_{l=1}^{M-1} \Xi(z_l)\mathds 1_{[z_{l},z_{l+1}]}$ and $X=\Xi(z_M)$, 
we identify ${\mathbf J}$ and $\mathbf I.$ Thus the proof of the proposition is complete.
\end{pr}

To complete the proof of  Theorem \ref{theo1}(1),
we now need to show that the LDP is also true for the uniform topology.
From Proposition \ref{DG} and Lemma \ref{lemtight},
and as the topology of uniform convergence is finer than the topology of 
pointwise convergence,  we can apply \cite[Corollary 4.2.6]{DZ} and get that 
the law of $\left((K^{n}(z))_{z\in \K}, C^{n}\right)$
as an element of  $\mathcal C(\K, {\Hr})\times {\Hr}$ equipped with the uniform topology 
satisfies a LDP in the scale $n$ with good rate function $\bf J.$ \\

\subsection{Properties
 of the rate function} \mbox{}\\
\label{subident}

To finish the proof of Theorem \ref{theo1}(1), the last thing to check 
is that $\mathbf I(K(\cdot),C)$ is infinite
 whenever $K$ is not Lipschitz continuous. This is the object of this subsection (see Lemma \ref{proprf}.(6)),
together with providing
further information on the functions $(K,C)$ with finite $\mathbf I$ that will be useful in the sequel.\\

We will
 consider the operator norm, given, for $H \in {\Hr}$, by $\|H\|_\infty= \sup  \langle u, Hu\rangle$,
where the supremum is taken over vectors $u\in \C^r$ with norm one. We also use the usual order on Hermitian matrices,
i.e. $H_1 \le H_2$ if and only if $H_2-H_1$ is positive semi-definite (respectively $H_1 <H_2$ if $H_2-H_1$
is positive definite).\\
 We recall that $\Lambda$ was defined in \eqref{1.02.10.def_La}.

\begin{lem}\label{proprf}
\begin{enumerate}
\item $H\mapsto\Lambda(H)$ is increasing, $\Lambda(-H) \le 0$ if $H \ge 0$.
\item If we denote by  $(C^*)_{ij}={\mathbb E}[\overline {g_i}g_j]$. Then, for any $H \in {\Hr},$
$$\Lambda(H)\ge \Tr(H C^*).$$
\end{enumerate}

If we assume moreover that $G$ satisfies the first part of Assumption \ref{hyponG} (existence of some exponential moments),
we have the following properties.\\

\begin{enumerate}
\setcounter{enumi}{2}
 \item There exists $\gamma>0$ so that 
$$B:=\sup_{H:\|H\|_\infty\le \gamma}
\Lambda(H)<\infty\,.$$
\item  If ${\bf I}(K(\cdot),C)$  is finite, $C\ge 0$ and $K(z)\ge 0$, for   any $z \in \K.$
Moreover, for all $L$, there exists a finite constant $M_L$
so that on
$\{{\bf I}\le L\}$,  we have 
$$  \sup_{z\in {\mathcal K}}\|K(z)\|_\infty\le M_L, \quad  
\|C\|_\infty
\le M_L.$$
\item  If ${\bf I}(K(\cdot),C)$ or ${\bf J}(K(\cdot),C)$ are finite, 
then $z\ra K(z)$  is non increasing. 
\item For all $L$, there exists a finite constant $M_L$
so that on
$\{{\bf I}\le L\}$,  we have 
$$ \textrm{and for all }
z_1,z_2\in\K,
\|K(z_2)-K(z_1)\|_\infty \le M_L |z_1-z_2|\,.$$
In particular, $K'$ exists almost surely and is bounded by
$M_L$. 

\end{enumerate}

If we  assume now that $G$ satisfies both parts of Assumption \ref{hyponG} (the law of $G$ does not put mass on hyperplanes),
we then have the following additionnal properties.

\begin{enumerate}
\setcounter{enumi}{6}
\item
For all non null positive semi-definite $H\in \Hr$,  \be\label{La_inf_inf.1.02.10}\lim_{t\to+\infty} \La(-tH)=-\infty.\ee  
\item If ${\bf I}(K(\cdot),C)$ is finite, then $C >0$ and  $K(z) >0$ for   any $z \in \K.$ Moreover, for almost any $z \in \K$ and for any non zero vector $e$,  there is no  interval with non-empty interior on which  the function 
 $\langle e, K'(.) e\rangle$ vanishes everywhere.
\end{enumerate}
\end{lem}

\begin{pr}
\begin{enumerate}
\item The first  point is just based on the fact that almost surely,  $\Tr(HZ) \ge 0$ if $H\ge 0$. 
\item  The second point follows from Jensen's
inequality.
\item The third point is  due to the fact that
$\Tr(HZ)\le \|H\|_\infty\sum_{i=1}^r |g_i|^2$ so that 
by H\"older's inequality,
$$\Lambda(H)\le \log \E[e^{\|H\|_\infty\sum_{i=1}^r |g_i|^2}]
\le  \frac{1}{r}
 \sum_{i=1}^r \log \E[e^{\|H\|_\infty r |g_i|^2}]$$
which is finite by Assumption \ref{hyponG} if $\|H\|_\infty r\le \alpha.$\\
\item To prove the fourth point let $(C,K)\in \{{\mathbf I}\le L\}$.
We first show that $C\ge 0$. 
We  take $P,X\equiv 0$ to get
$$\sup_{Y\in {\Hr}}\{\Tr(C Y)- \Lambda(Y)\}\le {\bf I}(K,C)\le L\,.$$
Suppose now that there exists some vector $u\in \mathbb C^r$ such that  $\langle u, C u\rangle=\alpha<0 ,$ and define, for any $t>0,$ $Y_t=-t\, uu^*$.
Then $\Lambda(Y_t)\le 0$ by the first point
and $\Tr(C Y_t) = -\alpha t$ so that for all $t>0$, 
$$-\alpha t \le \Tr(C Y_t)- \Lambda(Y_t)\le L.$$
Letting $t$ going to infinity gives a contradiction.
The same proof holds for $K(z)$ by taking $P(z)=-\one_{z\ge z_0} X$
and $X=-t\, uu^*$ if $\langle u, K(z_0) u\rangle=\alpha<0 .$
We finally bound $K$ and $C$.
With $\gamma$ and $B$ introduced in the third point, we define
 $Y=\pm \gamma uu^*$ and take $P,X\equiv 0.$ We get
$$\gamma |\langle u, C u\rangle  |\le B+L$$
for all vector $u$ with norm one, that is $\|C\|_\infty \le
\gamma^{-1}( L+B).$ Similar considerations hold for the bound over $\|K(z)\|_\infty.$\\

\item We next prove that $z\ra K(z)$ is non increasing when the entropy
is finite.
 Let us prove that for any $z_1, z_2\in \K$ 
\st $z_1<z_2$, $K(z_2)\le K(z_1)$ (dividing by $z_2-z_1$ will then give the fact that $K'$ is negative semi-definite where it is defined).
So let us fix $z_1, z_2\in \K$ 
\st $z_1<z_2$. Let us fix $u\in \C^r\bck\{0\}$. For all  real number  $t\ge 0$,  we have, for $P_t(z):=t \one_{[z_1, z_2]}(z)uu^*$ and $X=Y=0$, $$I(K(\cdot),C)\ge t u^*(K(z_2)-K(z_1))u -\Gamma(P_t,0,0).$$ 
Note that \bes \Gamma(P_t,0,0)=\int \La \lf(-t\int_{z_1}^{z_2}\f{\ud z}{(z-x)^2}u^*u\ri)\ud \mu(x)\le 0\ees
by (1) of this lemma. Thus for all $t>0$, $$I(K(\cdot),C)\ge t u^*(K(z_2)-K(z_1))u.$$ It follows that $u^*(K(z_2)-K(z_1))u$ is non positive
by letting $t$ going to infinity, which completes the proof of this point.

\item 
Take  $P= -(z_2-z_1)^{-1} 1_{[z_1,z_2]}uu^*$,
$Y=-uu^* \max_{x\in \mbox{supp}(\mu)} \int (z-x)^{-2} (z_2-z_1)^{-1} 1_{[z_1,z_2]}(z) dz$ and $X=0$ to get, since then $\tilde\Gamma(P,Y,X)\le 0$ by the first point,
$$\langle u, -\Tr((K(z_2)-K(z_1))(z_2-z_1)^{-1}) u\rangle \le L +r\e^{-2}\|C\|_\infty$$ 
where we used that $Y$ is bounded by $\e^{-2}$. This 
 provides the expected bound by the fourth point.

\item   Consider $\eta>0$ and a non vanishing
 orthogonal projector $p\in \Hr$ \st $H\ge \eta p$.   For all $t>0$, we have 
$$0\le \E[e^{-t\Tr( HZ)}]\le  \E[e^{-t\eta \Tr (pZ)}]= \E[e^{-t\eta \Tr (pGG^*)}]= \E[e^{-t\eta G^*pG}].$$ Since, by dominated convergence, 
 $$\lim_{t\to+\infty} \E[e^{-t\eta G^*pG}]= \Pro\{G^*pG=0\}=\Pro\{G\in \ker p\}=0$$(where we used   Assumption \ref{hyponG} in the last equality), we have $$\lim_{t\to+\infty} \La(-tH)=\lim_{t\to+\infty} \log\E[e^{-t\Tr( HZ)}]=
 -\infty.$$

\item We already proved that $K$ is
non increasing and almost surely differentiable, so that $K'\le 0$ almost surely. Moreover, if $u$ is a fixed vector and $\langle u, K'(\cdot) u\rangle$
vanishes on an interval $[z_1,z_2]$ with $z_1<z_2,$ taking 
$P_t=t\one_{[z_1,z_2]}(z)uu^*$, and $X=Y=0$, yields
$$I(K(\cdot),C)\ge  -\int \Lambda\left(-t\int_{z_1}^{z_2} \frac{dz}{(z-x)^2} 
uu^*\right)d\mu(x)$$
which goes to infinity as $t$ goes to infinity by the previous consideration.
Thus, this is not possible.  As we have already seen that $K(a)\ge 0$
for all $a\in \K$, we see that $K(a')>0$ for $a'<a$ unless 
there exists $e$ so that $\langle e, (a-a')^{-1}(K(a)-K(a')) e\rangle$
vanishes, which is impossible by the above.
\\

\end{enumerate}
 \end{pr}

\subsection{Study of the minimisers of $\mathbf I$}\mbox{}\\

We characterise the minima of ${\mathbf I}$ as
follows~:
\begin{lem}\label{lemmin}
For any compact set $\K$ of $(b,\infty)$,
the unique minimizer of ${\mathbf I}$ 
on  $\mathcal C(\K, {\Hr})\times {\Hr}$ 
is the pair $(K^*, C^*)$   given, for $1\le i , j \le r$, by
$$(K^*(z))_{ij}=\int\frac{(C^*)_{ij}}{z-\lambda}\ud\mu(\lambda), \,\,\,  \textrm{ for } \,\,\, z\in\K
\quad \textrm{ and } \quad(C^*)_{ij}=\mathbb E[\overline{g_i} g_j] .$$

\end{lem}

\begin{pr}
${\mathbf I}$ vanishes at its minimisers (as a good rate
function) and therefore a minimizer $(K,C)$
satisfies for all $P,X,Y$,
\begin{equation}\label{toto}
\Tr\left(\int K' (z)P(z)\ud z+K(z^*)X+C Y\right)
\le \Gamma(P,X,Y)\,.
\end{equation}
Now, for any fixed $(P,X,Y),$ there exists $\e_0 >0$ such that for any   $0<\e <\e_0,$
for any $x$ in the support of $\mu$ we have
$$ \e\left\| -\int \frac{1}{(z-x)^2}P(z)\ud z+ \frac{1}{z^*-x}X +Y\right\|_\infty < \alpha,$$
with $\alpha$ given by Assumption \ref{hyponG}. Therefore,  there exists 
a constant $L$ such that for any $x$ in the support of $\mu$
\begin{multline*}
 \left|\E\left(e^{\e\Tr\left(-\int \frac{1}{(z-x)^2}P(z)\ud z + \frac{1}{z^*-x}X+Y\right)Z}\right) \right.\\
\left.- \E \left(1 +\e \Tr \left(\left(-\int \frac{1}{(z-x)^2}P(z)\ud z  + \frac{1}{z^*-x}X+Y\right)Z\right)\right)\right| \le \e^2 L,
\end{multline*}
so that
$$\Gamma(\e P,\e X, \e Y)=\e \Tr\left(\int (K^*)' (z)P(z)\ud z+K^*(z^*)X + C^*
Y\right)+O(\e^2).
$$
As a consequence, for any minimizer $(K,C),$ we find after replacing $(P,X,Y)$
by $\e (P,X,Y)$, using \eqref{toto}
and letting $\e$ going to zero, that
$$\Tr\left(\int K' (z)P(z)\ud z+K(z^*)X+ C Y\right)
\le 
\Tr\left(\int (K^*)' (z)P(z)\ud z+K^*(z^*)X+C^* Y\right)\,.$$
Changing $(P,X,Y)$ in $-(P,X,Y)$
gives the  equality. This implies that
$$C=C^*,\qquad  K'=(K^*)'\quad \mbox{ a.s. and } \quad K(z^*)= K^*(z^*)$$
and therefore $(K,C)=(K^*,C^*).$
\end{pr}

\section{Large deviations for the largest eigenvalues in the  case without outliers}
\label{SectionLDPvap}

\noindent
We again assume throughout this section that Assumptions \ref{spacing}, \ref{hyponG} and \ref{without} hold.

\subsection{Statement of the main result}\mbox{}\\

 For any $\e >0, $ we define the compact set 
$\K_\e:= [b+\e,\e^{-1}].$ Let $s:= \operatorname{sign}\left(\prod_{i=1}^r \theta_i\right)=(-1)^{r-m}$.\\
For $x\in \R$, we set $\R^p_\downarrow(x) = \{(\alpha_1, \ldots, \alpha_p) \in \R^p / \alpha_1 \ge \ldots \ge \alpha_p \ge x\}.$\\
We also denote by $\omega(g) := \sup_{x \neq y} \frac{|g(x)-g(y)|}{|x-y|} \in [0,\infty]$ the  Lipschitz constant of a function $g.$ 
For any $\e, \gamma > 0,$ and $\alpha \in \R^p_\downarrow (b+\e),$ we put 
\begin{multline*}
S_{ \alpha, \gamma}^\e:=\left\{f\in \mathcal C(\K_\e, \R):  \exists
 g \in \mathcal C(\K_\e,\R) \textrm{ with }   \gamma \le g\le  \frac 1\gamma, \omega(g) \le  \frac 1 \gamma  \right.\\
\left.\textrm{ and }
f(z)=s.g(z) \prod_{i=1}^p (z-\alpha_i)\right\}
\end{multline*}
Note that in the latter product, the $\alpha_i$'s appear with multiplicity.
$S_{ \emptyset, \gamma}^\e$ will denote the set of functions as above but with
no zeroes on $\K_\e$. 
We have the following theorem.
\begin{Th}\label{main}
Under Assumptions \ref{spacing}, \ref{hyponG} and \ref{without},
 the law of the $m$ largest eigenvalues $(\wtl_1^{n},\ldots, \wtl_m^{n})$ of $\wtX$  satisfies
a large deviation principle in $\R^m$ in the scale
$n$ and  with good rate function $L,$ defined as follows. For $\mathbf \alpha = (\alpha_1, \ldots, \alpha_m)\in \R^m,$ 
we take $\alpha_{m+1} =b$ and  
\[
L( \alpha ) := \begin{cases}

\lim_{\e \downarrow 0}\inf_{\cup_{\gamma >0} S_{(\alpha_1, \ldots, \alpha_{m-k}), \gamma}^\e} J_{\K_\e} &
\textrm{if } \alpha \in \R^m_\downarrow(b), \alpha_{m-k+1} = b \textrm{ and } \alpha_{m-k} >b   \\
& \qquad\textrm{ for some }
k\in \{0, \ldots, m-1\} ,\\
\lim_{\e \downarrow 0}
\inf_{\cup_{\gamma >0} S_{\emptyset, \gamma}^\e} J_{\K_\e} &
\textrm{if } \alpha_1=\alpha_2=\cdots\alpha_m=b \\
+\infty & \textrm{otherwise.}
\end{cases}\]
\end{Th}

\begin{rmq} \label{increasing} The function $L$ is well defined. Indeed, one 
 can easily notice that for all $ \alpha \in \R^m_\downarrow(b) 
$ \st for some $k\in \{0,\ldots, m\}$,  $\alpha_{m-k+1} = b$  and
$ \alpha_{m-k} >b$, the map $$\eps\longmapsto \inf\{J_{\K_\e} (f)\ste f\in {\cup_{\gamma >0} S_{(\alpha_1, \ldots, \alpha_{m-k}), \gamma}^\e}\}$$ is increasing,
so that its limits as $\e$ decreases to zero exists.\end{rmq}
\begin{rmq} \label{rmqweyl}
Note that $ J_{\K_\e}(f)$ is infinite if $f$ has more
than $r$ zeroes  greater than  $b$. Indeed, by definition, if $ J_{\K_\e}(f)$ is finite,
$$f(z)=  P_{\Theta,r}(K(z),C)=c\det(A-K(z))$$
with a non-vanishing constant $c$ and a self-adjoint matrix $A$
  with eigenvalues $(\theta_1^{-1},
\ldots,\theta_r^{-1})$ and a function  $K$  with
values in the set of  $r\times r$ positive self-adjoint matrices so that $K'\le 0$ by
Lemma \ref{proprf}.  We may assume without loss of
generality that $f$ vanishes at a point $x>b$, since otherwise we
are done, so that
there exists a non zero  $e\in\mathbb C^r$ so that $K(x)e= Ae$.
There is  at most one $x$ at which   $K(x)e= Ae$; otherwise, 
 $\langle e, K'(\cdot) e\rangle$ would vanish on a non trivial interval 
which is impossible by (7) of Lemma \ref{proprf}.
Moreover, if we let $P$ be the orthogonal projection
onto the orthocomplement of $e$, the function $H(z)=\det( (1-P)(A-K(z))(1-P))
 \det(PAP-PK(z)P)$    vanishes at $x$
and at the zeroes of $ \det(PAP-PK(z)P)$. But 
 $PAP$ and $PK(z)P$
have the same properties as $A$ and $K(z)$ except they
have one dimension less. Thus, we can proceed by induction
and see that $f$ can vanish at at most $r$ points. 

\end{rmq}

The minimisers are described by the following result.
\begin{Th}\label{min}
If we define on $(b,\infty)$
$$H(z)=P_{\Theta,r}(K^*(z), C^*)$$
 where $(K^*,C^*)$ are given in Lemma
\ref{lemmin} and $P_{\Theta,r}$ is defined in Proposition \ref{zeroesscheme}, there exists $k\in \{0, \ldots, m\}$ such that $H$
has $m-k$ zeroes $(\lambda_1^*, \ldots, \lambda_{m-k}^*) $ (counted with multiplicity). 
The unique point of $\R^m$ on which  
$L$ vanishes is  $(\lambda_1^*, \ldots, \lambda_{m-k}^*,b ,\ldots,b)$ and consequently
$(\wtl_1^{n},\ldots, \wtl_m^{n})$ converges almost surely to this point as $n$ grows to infinity.
\end{Th}

\begin{rmq} \label{rmqmin}
 In the case when $(g_1,\ldots, g_r)$ are independent centered variables
with variance one, one can check that 
$C^* = I_r,$ $K^*(z)= \int \frac{1}{z-x} \ud\mu(x).I_r$ and $$H(z) = \prod_{i=1}^r \left(\frac{1}{\theta_i}-\int \frac{1}{z-x} \ud\mu(x) \right)$$
so that we recover \cite[Theorem 2.1]{benaych-rao.09} or \cite[Theorem 1.3]{benaych-guionnet-maida_fluctu}.
\end{rmq}

\subsection{Preliminary remarks and strategy of the proof}\mbox{}\\

Let us first notice that at most $m$ eigenvalues
of $\wtX$ can deviate from the bulk since 
 by 
Weyl's interlacing inequalities (see e.g. \cite[Section 4.3]{Horn})
$$ \wtl^{n}_{m+1} \le  \la^{n}_1  
 , $$ which converges to  $b$ as  $n$ goes to infinity.

Secondly, let us state the following lemma.
\begin{lem}\label{lemexpt}
The law of the sequence $(\wtl_1^{n},\ldots, \wtl_m^{n})$ of the $m$ largest eigenvalues of $\wtX$ is exponentially tight in the scale $n$.
 \end{lem}

\begin{pr} Let us define $R_n:=\wtX-X_n$ and  denote by $\|R_n\|_\infty$ the operator norm of the  perturbation matrix $R_n$. Note that for all $k$,
$$ \la^{n}_{k}-\|R_n\|_\infty \le \wtl_k^{n} \leq  \la^{n}_k +\|R_n\|_\infty. $$
Since for any fixed $k$, the non random sequence $\la_k^n$ converges to $b$ as $n$ tends to infinity, it suffices to prove that 
\be\label{exp_tight_30.06.10}\limsup_{L\to\infty}\limsup_{n\to\infty}\ff{n}\log\Pro(\|R_n\|_\infty\ge L)=-\infty.\ee
For the orthonormalized perturbation model, since $\|R_n\|_\infty=\max \{\theta_1, -\theta_r\}$, \eqref{exp_tight_30.06.10} is clear.
In the i.i.d. perturbation model, we have, for $\theta:=\max_{1\le i\le r}|\theta_i|$, 
$$
\|R_n\|_\infty=\sup_{\|v\|_2=1}|\langle v, R_n v\rangle|
\le \ff{n} \sum_{i=1}^r\theta\|G_i^n\|_2^2=\f{\theta}{n}\sum_{k=1}^n \sum_{i=1}^r|g_i(k)|^2.$$
It implies, by Tchebychev's inequality, that
$$\Pro\left(\|R_n\|_\infty\ge L\right)\le e^{-\f{n\alpha L}{\theta}}\E\left[\exp\left(\alpha\sum_{k=1}^n \sum_{i=1}^r|g_i(k)|^2\right)\right]=e^{-\f{n\alpha L}{\theta}}\E\left[\exp\left(\alpha \sum_{i=1}^r|g_i(1)|^2\right)\right]^n$$
which allows to conclude by  Assumption \ref{hyponG}.
\end{pr}

As the  law of  $(\wtl_1^{n},\ldots, \wtl_m^{n})$
is exponentially tight, 
  the proof of Theorem \ref{main}   reduces to  establishing a weak LDP.
In virtue of   \cite[Theorem 4.1.11]{DZ} (see also \cite[Corollary D.6]{alice-greg-ofer}), this weak LDP (and the fact that
$L$ is a rate function) will be a direct consequence 
of Equation \eqref{outside} and  Lemma \ref{near} below. The fact that $L$ is a good rate function is then implied 
by exponential tightness \cite[Lemma 1.2.18]{DZ}. \\

\subsection{The structure of $H^{n}$}\mbox{}\\

From Proposition \ref{zeroesscheme}, we know that the $\wtl^{n}_i$'s are essentially the zeroes of $H^{n}.$
However, $H^{n}$ could a priori have other zeroes than these eigenvalues or
take arbitrary  small values.
To control this point, we need to understand better the structure of $H^{n}$.
Let
\begin{multline*}
C_{k,\gamma}^\e:=\bigg\{f\in \mathcal C(\K_\e,\R): \exists p\textrm{ polynomial of
degree $k$ with } k \textrm{ roots in } \K_\e  \\
\textrm{ and dominant coefficient } 1,
 g \in \mathcal C(\K_\e,\R) \textrm{ with }  \gamma \le g\le  \frac 1\gamma, \omega(g) \le  \frac 1 \gamma  \textrm{ and }
f(z)=s.g(z)p(z)\bigg\}
\end{multline*}
and 
$$C_{\gamma}^\e = \bigcup_{0 \le k \le m }C_{k,\gamma}^\e.$$  
We intend to show the following fact
\begin{lem}\label{lemonH}
For any $\e >0$ small enough, there exists a positive integer $n_0(\e)$, $L(\e) >0$ and
a sequence of random functions $(g_n)$ such that
for any $z \in \K_\e$ and $n \ge n_0(\e),$
 \begin{equation*}\label{form}
H^{n}(z)=\begin{cases}s\prod_{i=1}^r \|q_i^{n}W_i^{n}\|_2^2 \,\,
g_n(z)\,\,
\prod_{i=1}^m (z-\wtl_i^{n})&\textrm{in the orthonormalized perturbation model,}\\
&\\
s \,\,
g_n(z)\,\,
\prod_{i=1}^m (z-\wtl_i^{n})&\textrm{in the i.i.d. perturbation model,}\end{cases}
\end{equation*}
with
\be \label{encadrgn}
  s=(-1)^{r-m}, \qquad L(\e) \le g_n \le \frac{1}{L(\e)} \qquad\textrm{ and }\qquad \omega(g_n) \le \frac{1}{L(\e)}.\ee
In particular, for any $\e >0,$
\be \label{Htight2}
 \limsup_{\gamma\downarrow 0}
\limsup_{n\ra\infty}\frac{1}{n}\log \Pro\left( 
\left(\prod_{i=1}^m (z-\wtl_i^{n})^{-1}
H^{n}(z)\right)_{z \in \K_\e} \in (C_{0,\gamma}^\e)^c\right)= -\infty.
\ee 
and

\be \label{Htight}
 \limsup_{\gamma\downarrow 0}
\limsup_{n\ra\infty}\frac{1}{n}\log \Pro\left( 
H^{n}\in (C^\e_\gamma)^c\right)= -\infty.
\ee  

\end{lem}

\begin{pr}
Let us define the random sequence $$c_n:= \begin{cases}s\prod_{i=1}^r \|q_i^{n}W_i^{n}\|_2^2  &\textrm{in the orthonormalized perturbation model,}\\
&\\
s&\textrm{in the i.i.d. perturbation model.}\end{cases}$$ Going back to the proof of Proposition \ref{zeroesscheme}, one can easily see that, for any $z \notin \{\la_1^{n}, \ldots, \la_n^{n}\},$
$$
H^{n}(z)=c_n\prod_{i=1}^r \theta_i^{-1} \det(zI_n-X_n)^{-1}\det\left(zI_n-X_n-\sum_{i=1}^r\theta_i U_i^{n}(U^{n}_i)^*\right).$$
We can rewrite the above as
$
H^{n}(z)=c_n
g_n(z)\,\,
\prod_{i=1}^m (z-\wtl_i^{n})$
with 
$$
 g_n(z):= \prod_{i=1}^r |\theta_i|^{-1}
\frac{1}{ \prod_{i=1}^m (z-\lambda_i^{n})} 
\prod_{i=m+1}^{n} \left(1+
\frac{\lambda_i^{n} -\wtl_i^{n}}{ z-\lambda_i^{n}}
\right) $$
Now, for $\e >0$  fixed, we shall bound $g_n$ and its Lipschitz constant on $\K_\e$.

As $\K_\e$ is compact and the $\lambda_i^{n}$
belong to a fixed compact, for $\e >0$ small enough, 
 for any $i$ and $z\in \K_\e$ we have
 $ z-\lambda_i^{n} \le \frac 2 \e$  and $|\lambda_i^{n}| \le \frac 2 \e$ so that
\be\label{cont}
 0\le \sum_{i=m+1}^{n} ( \lambda_{i-m}^{n}- \lambda_{i}^{n})
=\sum_{i=1}^m \lambda_i^{n} -\sum_{i=n-m}^n \lambda_i^{n} \le 2m \frac 2 \e.\ee
 We choose $n_0(\e)$ such that  for $n \ge n_0(\e)$ and any $i$ and $z\in \K_\e$
we have $\frac \e 2 \le z-\lambda_i^{n}$ so that as $ z-\lambda_i^{n} \le \frac 2 \e$,
 $$  0 \le \frac{\la_{i-m}^{n}- \la_{i}^{n}}{ z-\lambda_i^{n}} =
1- \frac{z-\la_{i-m}^{n}}{ z-\lambda_i^{n}} 
\le
1-\frac{\e^2}{4}.$$

Now, using Weyl's interlacing properties, we have
 for any $i \ge m+1,$
$$\wtl_i^{n} \le   \lambda_{i-m}^{n},\mbox{ 
so that
} \lambda_i^{n} -\wtl_i^{n} \ge  -( \lambda_{i-m}^{n}- \lambda_{i}^{n}).$$

For $0 \le x \le 1-\frac{\e^2}{4},$ $\log(1-x) \ge - \frac{4}{\e^2} x,$ so that
we finally get by \eqref{cont},
$$ g_n(z)
\ge \prod_{i=1}^r |\theta_i|^{-1} \left(\frac{\e}{2}\right)^{m} 
 e^{ \frac{4}{\e^2}\sum_{i=m+1}^{n} 
 \frac{\lambda_i^{n} -\wtl_i^{n}}{ z-\lambda_i^{n}}}
\ge
\prod_{i=1}^r |\theta_i|^{-1} \left(\frac{\e}{2}\right)^{m} 
e^{- 2m \left(\frac{4}{ \e^2}\right)^2}.$$
By very similar arguments (using $\log(1+x ) \le x$ for $x\ge 0$), one can also check that for any $n \ge n_0(\e),$
$$ g_n(z)
\le \prod_{i=1}^r |\theta_i|^{-1} \left(\frac{2}{\e}\right)^{m} e^{\frac{4}{\e^2}2m}. $$
The proof of the
 uniform equicontinuity of $g_n$ on $\K_\e$ is left to the reader
as the arguments are very similar since $z\ra (z-\lambda_i^n)^{-1}$
is uniformly continuous on $\K_\e$ for $n$ large enough. 

To prove \eqref{Htight} and \eqref{Htight2}, it is therefore
sufficient 
to prove that, with probability greater than $1- e^{-cn}$ for somme $c>0,$ we have that $c_n$ and $c_n^{-1}$ are bounded  which is a direct consequence of Lemma \ref{norm},
and that, $\tilde\lambda_n^1 \le \e^{-1}$
for small $\e$, which is
proved in Lemma \ref{lemexpt}.
\end{pr}

The main application of the previous Lemma will be the following
continuity properties of the zeroes of functions in $C^\e_\gamma$.
\begin{lem}\label{contz}
Let $\e>0$ be fixed, $\gamma>0$ small enough,
 and $k\in\mathbb N$ be fixed. Let $\alpha_1^0\ge \cdots\ge\alpha^0_k\in \K_\e $ and $f_0
(z)=h_0(z)\prod_{i=1}^k (z-\alpha_i^0)\in C^\e_{k,\gamma}$, 
 be given.
Then, for all $\delta>0$ there exists $\delta'>0$ so that
\begin{multline*}
\left\{
f\in C^\e_\gamma:\sup_{x\in\K_\e}|f(x)-f_0(x)|<\delta'
 \right\}\\
\subset \left\{ z\mapsto h(z)\prod_{i=1}^m (z-\alpha_i):
h\in C^\e_{0,\gamma}, \max_{1\le i\le k}|\alpha_i-\alpha_i^0|\le \delta,
\max_{i> k}\alpha_i\le b+2\e\right\}
\end{multline*}
\end{lem}
\begin{pr}
This amounts to show that if $f^n\in C^\e_\gamma$
is a sequence converging (for the uniform topology on
$\K_\e$) to $f\in  C^\e_{k,\gamma}$, $m-k$
zeroes of the functions $f_n$ will be below $b+2\e$
and the others will converge to the zeroes of $f$. 
Indeed, if we take a sequence
$f^n\in C^\e_\gamma$, we can always denote it 
$f^n(z)=h^n(z)\prod_{i=1}^m (z-\alpha_i^n)$ (with possibly 
some $\alpha_i^n\in (b+\e/4, b+3\e/4)$  if $f^n\in C_{k,\gamma}^\e$
with $k<m$), as this amounts
at the worst to change $h$ and
 take $\gamma$ smaller.  Then, the crucial point is that  $h^n$ is tight by Arzela-Ascoli theorem so that 
 we can consider a converging subsequence. As the
$\alpha^n_i$ belong to $[b,1/\e]$, we can also consider converging subsequences.
Thus, $f^n$ converges along subsequences to a function $\tilde f$ with $\tilde f(z)=h(z)\prod_{i=1}^m (z-\alpha_i)$  on $\K_\e$ with $\alpha_i\in [b,1/\e]$. But then we must have
$f=\tilde f$ which allows in particular to identify $k$ limit points with the zeroes of $f$, the others being below $b+2\e$.
\end{pr}

\subsection{Core of the proof}\mbox{}\\

\label{core}

First, from what we said in the preliminary remarks and the fact that the $\wtl^{n}_i$ are decreasing,
we obviously have that  if $\alpha \notin  \R^m_\downarrow(b)$, one has
\begin{multline}\label{outside}
 \limsup_{\delta\downarrow 0}
\limsup_{n\ra\infty}\frac{1}{n}\log \Pro\left( \bigcap_{1\le i\le m}\{|\wtl^{n}_i-\alpha_{i}|\le\delta\}\right)\\
=\liminf_{\delta\downarrow 0}
\liminf_{n\ra\infty}\frac{1}{n}\log \Pro\left( \bigcap_{1\le i\le m}
\{|\wtl^{n}_i-\alpha_{i}|<\delta\}\right)=-\infty. 
\end{multline}

The weak LDP will then be a direct consequence of the following lemma, with $k$ the numbers of eigenvalues going to $b,$

\begin{lem} \label{near}
Let $\alpha \in \R^m_\downarrow$ and $k$ between $0$ and $ m$ such that $\alpha_{m-k+1} = \ldots = \alpha_{m}=b$ and $\alpha_{m-k}>b$ if $k<m$.
We have
\begin{multline*}
\lim_{\e \downarrow 0}\limsup_{\delta \downarrow 0}
\limsup_{n\ra\infty}\frac{1}{n}\log \Pro\left( \bigcap_{1\le i\le m-k}\{|\wtl^{n}_i-\alpha_{i}|\le\delta\}\bigcap_{m-k+1\le i\le m}\{\wtl^{n}_i\le b+\e \}\right)
\\
=\lim_{\e \downarrow 0}\liminf_{\delta \downarrow 0}
\liminf_{n\ra\infty}\frac{1}{n}\log \Pro\left(  \bigcap_{1\le i\le m-k}\{|\wtl^{n}_i-\alpha_{i}|\le\delta\}\bigcap_{m-k+1\le i\le m}\{\wtl^{n}_i\le b+\e \}\right)
=-L(\alpha),
\end{multline*}
with the obvious convention that $\bigcap_{m-k+1\le i\le m}\{\wtl^{n}_i\le b+\e \}= \Omega$ if $k=0.$
\end{lem}

\begin{pr}
Let $\delta$ and $\e$ be  positive small enough constants so that 
 $\alpha_{m-k}-\delta \ge b+2\e.$ In particular,
$\cap_{i=1}^{m-k} [\alpha_i-\delta,\alpha_i+\delta] \subset \K_\e.$
On the set $\bigcap_{1\le i\le m-k}\{|\wtl^{n}_i-\alpha_{i}|\le\delta\}\bigcap_{m-k+1\le i\le m}\{\wtl^{n}_i\le b+\e \}$,
 for all $i \le m-k,$ $\wtl^{n}_i$ is in $\K_\e.$ On the other hand,
 for $n$ large enough,  $\{\la_1^{n}, \ldots, \la_n^{n}\}\cap \K_\e =   \emptyset.$
Therefore,  $\wtl^{n}_i \notin \{\la_1^{n}, \ldots, \la_n^{n}\}$
for $i\in\{1,\ldots,m-k\}$ and,
by Proposition \ref{zeroesscheme}, is a zero of $H^{n}.$ 

Let us next prove the large deviation upper bound
and fix $ \alpha_1\ge\alpha_2\ge\cdots\ge\alpha_{m-k}>b$.   A function $f\in  C_{\gamma,k}^\e$ which vanishes
within a distance $\delta$ of $(\alpha_i)_{1\le i\le m-k}$ with $\delta<\alpha_{m-k}-b$
belongs to the set 
\begin{multline} \label{Bkde}
 B_{\alpha,\gamma,\delta}^\e: = \Bigg\{f\in \mathcal C(\K_\e,\R):  \exists
 g \in \mathcal C(\K_\e,\R) \textrm{ with }  \frac{1}{\gamma} \le g\le \gamma, \omega(g) \le \frac{1}{\gamma}\\
  \textrm{ and }
f(z)=s.g(z) \prod_{i=1}^{m-k} (z-\beta_i)
\textrm{ with } \forall i \le m-k, \beta_i \in [\alpha_i-\delta,\alpha_i+\delta] \Bigg\}
\end{multline}
Moreover, writing $H^n(z)=h^n(z)\prod_{i=1}^m (z-\tilde\lambda^n_i)$
by Lemma \ref{lemonH},
clearly $H^n$ belongs to  $ B_{\alpha,\gamma,\delta}^\e
$ as soon as  for some $\e'<\e$ and $\gamma'\cdot (\e')^m>\gamma$,
$h^n\in C_{0,\gamma'}^{\e'}$ and
$\bigcap_{1\le i\le m-k}\{|\wtl^{n}_i-\alpha_{i}|\le\delta\}\bigcap_{m-k+1\le i\le m}\{\wtl^{n}_i\le b+\e-\e' \}$ holds.
As a consequence, we can write 
\begin{multline*}
 \Pro\left(\bigcap_{1\le i\le m-k}\{|\wtl^{n}_i-\alpha_{i}|
\le\delta\}\bigcap_{m-k+1\le i\le m}\{\wtl^{n}_i\le b+\e-\e' \}\right)\\
\leq \Pro\left(H^{n} \in B_{\alpha,\gamma,\delta}^\e\right)+
\Pro\left(h^{n} \in (C_{0,\gamma'}^{\e'})^c\right).
\end{multline*}
Then, 
 by \cite[Lemma 1.2.15]{DZ},
\begin{eqnarray}
&& \limsup_{n\ra\infty}\frac{1}{n}\log \Pro\left( \bigcap_{1\le i\le m-k}\{|\wtl^{n}_i-\alpha_{i}|\le\delta\}
\bigcap_{m-k+1\le i\le m}\{\wtl^{n}_i\le b+\e -\e'\}\right)\nonumber\\ 
&&\quad\quad\le \max\left\{\limsup_{n\ra\infty}\frac{1}{n}\log \Pro\left(H^{n} \in B_{\alpha,\gamma,\delta}^\e\right); \limsup_{n\ra\infty}\frac{1}{n}\log
\Pro\left(h^{n} \in (C_{0,\gamma'}^{\e'})^c\right)
\right\},\label{polk}
\end{eqnarray}
 Moreover, $B_{\alpha,\gamma,\delta}^\e$ is a closed subset of $\mathcal C(\K_\e,\mathbb R).$ Indeed, if we take a converging 
sequence $f_n(z)=s g_n(z) \prod_{i=1}^{m-k}
(z-\beta_i^n)$, 
since the $\beta_i^n, n\ge 0$ belongs to compacts
and the $g_n, n\ge 0$ are tight by Ascoli-Arzela's theorem,
we can always assume up to extraction that $g_n$ and $\beta_i^n, 1\le i\le m-k$
converge so that the limit of $f_n$ belongs to $B_{\alpha,\gamma,\delta}^\e$.

Since $ J_{\K_\e}$ is a good rate function,
$(B_{\alpha,\gamma,\delta}^\e)_{\delta>0}$ is a nested family and
$\cap_{\delta>0} B_{\alpha,\gamma,\delta}^\e= S_{(\alpha_1, \ldots, \alpha_{m-k}), \gamma}^\e
$, Theorem \ref{theo1} gives with  \cite[Lemma 4.1.6]{DZ} that 
\be\label{leminterclosed}\limsup_{\delta\downarrow 0}\limsup_{n\ra\infty}
\frac{1}{n}\log \Pro\left(  H^{n} \in B_{\alpha,\gamma,\delta}^\e\right)\le -\inf_{ S_{(\alpha_1, \ldots, \alpha_{m-k}),\gamma}^\e } J_{\K_\e}.\ee
Taking $\gamma'=\gamma'_0$ small enough, \eqref{Htight2} and \eqref{polk} 
 give
for $\gamma/(\e')^m<\gamma_0'$,
\begin{multline*}
\limsup_{\delta\downarrow 0}
\limsup_{n\ra\infty}\frac{1}{n}\log \Pro\left( \bigcap_{1\le i\le m-k}\{|\wtl^{n}_i-\alpha_{i}|\le\delta\}
\bigcap_{m-k+1\le i\le m}\{\wtl^{n}_i\le b+\e -\e' \}\right)\\
\leq -\inf_{  S_{(\alpha_1, \ldots, \alpha_{m-k}),\gamma}^\e}  J_{\K_\e}
\leq -\inf_{ \cup_{\gamma>0} S_{(\alpha_1, \ldots, \alpha_{m-k}),\gamma}^\e}  J_{\K_\e}.
\end{multline*}
We can finally take $\gamma=0$ (nothing depends on it
anymore),  $\e-\e'$ going to zero,
as the left hand side  obviously decreases as $\e-\e'$
decreases to $0$ and, as we already mentioned it in Remark \ref{increasing},
 the right hand side increases as $\e$
decreases to $0.$\\

\noindent
We turn to the lower bound, which is a bit more delicate. Let us again consider  $\delta$ and $\e$ small enough so that
$\cap_{i=1}^{m-k} [\alpha_i-\delta,\alpha_i+\delta] \subset \K_\e$. As $  J_{\K_\e}$ is a good rate function and $S_{\alpha,\gamma}^\e$ is closed, for all $\gamma>0$, the infimum 
$\inf_{ S_{(\alpha_1, \ldots, \alpha_{m-k}),\gamma}^\e}  J_{\K_\e}$ is achieved,  say at $f_\gamma^{k,\e}.$
To complete the proof, we need the following lemma, 
based on the structure of $H_n$ and  whose  proof is a direct application of Lemma \ref{contz}.
\begin{lem} \label{zeromultiples}
Let $\e,\gamma$ be fixed and small enough. There exists $\delta_0$ such that for any $\delta \le \delta_0$, there exists 
$\delta^\prime >0$ \st for any $n,$  
$$\left\{H^n\in C^\e_\gamma\right\}\cap \left\{\sup_{x\in
\K_\e}| H^{n}(x)-f_\gamma^{k,\e}(x)|<\delta'\right\} \subset  \bigcap_{1\le i\le m-k}\{|\wtl^{n}_i-\alpha_{i}|\le\delta\}
\bigcap_{m-k+1\le i\le m}\{\wtl^{n}_i\le b+2\e \}. $$ 
\end{lem}

To prove the lower bound in  Theorem \ref{theo1}, we may assume without loss
of generality that $$J:=\lim_{\e\downarrow 0}
\inf_{ \cup_{\gamma >0} S_{(\alpha_1, \ldots, \alpha_{m-k}),\gamma}^\e}  J_{\K_\e}<\infty.$$
Let $\eta>0$ be fixed.
As $$\inf_{ \cup_{\gamma>0} S_{(\alpha_1, \ldots, \alpha_{m-k}),\gamma}^\e}  J_{\K_\e}
=\inf_{\gamma>0}\inf_{  S_{(\alpha_1, \ldots, \alpha_{m-k}),\gamma}^\e}  J_{\K_\e}=\inf_{\gamma>0}J_{\K_\e}(f_\gamma^{k,\e}), $$ 
 we can choose $\e, \gamma$ small enough so that 
$J_{\K_\e}(f_\gamma^{k,\e}) \le J+\eta$.
 By \eqref{Htight}, there exists $L(\gamma,\e)$ going to infinity as
$\gamma,\e$ go to zero so that for $n$ large enough,
$$ \Pro\left(H^n\in (C^\e_\gamma)^c\right)\le e^{-n L(\e,\gamma)}.$$
We choose $\gamma,\e$ small enough so that $L(\e,\gamma)>J+2\eta$.

Lemma \ref{zeromultiples} implies, that for $\delta \le \delta_0,$
for $\delta'$ small enough, $\eta>0$, 
for $n$ large enough,
\begin{multline*}\label{borne-inf.141209.15h21}
 \Pro\left( \bigcap_{1\le i\le m-k}\{|\wtl^{n}_i-\alpha_{i}|\le\delta\}
\bigcap_{m-k+1\le i\le m}\{\wtl^{n}_i\le b+2\e \}\right)\\ 
\ge \Pro\left(\sup_{z\in
\K_\e}| H^{n}(z)-f_\gamma^{k,\e}(z)|<\delta'\right)
- \Pro\left(H^n\in (C^\e_\gamma)^c\right)\ge \frac{1}{2}e^{-n(J+2\eta)}
\end{multline*}
the last inequality following from Theorem \ref{theo1}.(\ref{theoH}).
As $\eta$ can be chosen as small as we want, we conclude
by taking first $n$ going to infinity, and then $\delta,\e, \eta$
to zero.

\end{pr}

\subsection{Identification  of the minimizers}\mbox{}\\

We prove Theorem \ref{min}, 
which is straightforward. 
Since $L$ is a good rate function, it vanishes at its minimizers $(\lambda_1^*,\ldots,\lambda_m^*)\in \R^m_{\downarrow}(b).$
Putting $\lambda_{0}^*=b+1$,
we know that there exists $0\le k \le m$ such that $\lambda_{m-k}^*>b$
and $\lambda_{m-k+1}^*=b.$ From the definition of $L,$ for any $n$ large enough such that $b+\frac{1}{n}<\lambda_{m-k}^*, $ we can find 
a function $f_n$ defined on $\K_{\frac{1}{n}}$ vanishing at $(\lambda_1^*,\ldots,\lambda_{m-k}^*)$
such that $J_{\K_{\frac{1}{n}}}(f_n) \le \frac{1}{n}.$
 From the definition of $J_\K$ and the fourth and sixth point of Lemma \ref{proprf},
all the functions $f_n$  are in a compact set of $\mathcal C((b,\infty),\R)$
so that we can find a function $f$ vanishing at $(\lambda_1^*,\ldots,\lambda_{m-k}^*)$
so that $J_{\mathcal K_\e}(f)=0$
for all $\e>0$. But the latter  also implies that $f(z)=P_{\Theta,r}(K(z),C)$
with $(K,C)$ minimising ${\bf I}$, that is $(K,C)=(K^*,C^*)$
by Lemma \ref{lemmin}.

\section{Large deviations for the eigenvalues of Wishart
matrices}\label{wishartsec}

In this section, we study the i.i.d. perturbation model when $X_n=0.$
More precisely, we consider $G= (g_1, \ldots, g_r)$ satisfying Assumption \ref{hyponG},
$n\times r$ matrices $G_n$ 
whose rows are  i.i.d. copies of $G, $  a diagonal matrix 
$\Theta=\diag(\theta_1,\ldots,\theta_r)$ and we study 
 the large deviations of Wishart matrices  $W_n=\frac{1}{n} G_n \Theta G_n^*.$
This matrix has zero as an eiganvalue with muliplicity at least $n-r$ and we refer in the whole  section
to the $r$ eigenvalues of $W_n$ that can be non-zero as ``the eigenvalues of $W_n$''.
The large deviations
for the largest and smallest such eigenvalues were already studied in \cite{Fey}
in the case when $\Theta=(1,\ldots,1)$ and the $g_i$'s are i.i.d.

\begin{propo}\label{ldpwis}
Assume that $G$ satisfies Assumption \ref{hyponG}. Let $\Theta
=\diag(\theta_1,\theta_2,\ldots, \theta_r)$ be 
a diagonal matrix with positive entries. Then, the law of the eigenvalues
of $W_n$ satisfies a large deviation principle 
in the scale $n$ with rate function which is infinite unless $\alpha_1 \ge\cdots \ge\alpha_r \ge 0$ and in this case given by
$$L(\alpha_1,\ldots,\alpha_r)
=\inf\{ J(C)
: (\alpha_1,\ldots,\alpha_r) \mbox{ are the eigenvalues of } \Theta^{-\frac{1}{2}} C \Theta^{-\frac{1}{2}}\},$$
with
$$J(C)=\sup_{Y\in \Hr}\{\Tr(CY)-\log \E[e^{\langle G,Y G\rangle}]\}\,.$$
\end{propo}

Note that the previous proposition could also have been deduced directly from Cram\'er's theorem and the contraction principle. 

The Gaussian case allows an exact computation, given by the following 
\begin{cor}\label{ldpwisg}
Assume that $G = (g_1, \ldots, g_r)$ is a Gaussian vector with 
positive definite covariance matrix $R.$
 Let $\Theta
=\diag(\theta_1,\ldots, \theta_r)$ be 
a diagonal matrix with positive entries.
We denote by $0<r_1(\Theta) \le r_2(\Theta)\le\ldots \le r_r(\Theta)$ the eigenvalues
of the matrix $\Theta^{-1/2}R^{-1}\Theta^{-1/2}$ in increasing order.\\
 Then, the law of the  eigenvalues
of $W_n$ satisfies a large deviation principle 
in the scale $n$ with rate function which is infinite unless
$\alpha_1\ge\alpha_2\ge\ldots\ge\alpha_r >0$ and otherwise given by
$$L(\alpha_1,\ldots,\alpha_r)
=\frac{1}{2}\sum_{i=1}^r\lf( \alpha_i r_i(\Theta) -1 -\log (\alpha_i r_i(\Theta) )\ri)\,.$$
In the particular case when the entries are i.i.d. standard normal, the above rate function can be rewritten
 $$L(\alpha_1,\ldots,\alpha_r)
=\frac{1}{2}\sum_{i=1}^r\lf( \frac{\alpha_i} {\theta_i} -1 -\log \frac{\alpha_i} {\theta_i} \ri)\,.$$
\end{cor}

Now, by a straightforward use of the contraction principle, we can derive some results about the deviations of the largest eigenvalue. This problem
was addressed in particular in \cite{Fey}. The following corollary holds for the Gaussian case. 

\begin{cor}\label{lmawg}
 Under the assumptions of Corollary \ref{ldpwisg}, 
the
law of the largest eigenvalue satisfies a LDP with good rate function
$$
L_{\rm max}(x)=\begin{cases}
\frac{1}{2} (xr_1(\Theta) - 1 - \log ( xr_1(\Theta)))&\textrm{if } x \ge \frac{1}{ r_1(\Theta)},\\
\frac{1}{2} \sum_{i=1}^j (xr_i(\Theta) - 1 - \log ( xr_i(\Theta)))& \textrm{if } \frac{1}{ r_{j+1}(\Theta)} <  x \le  \frac{1}{ r_j(\Theta)}, 
\end{cases}$$
with the convention that $r_{r+1}(\Theta) = \infty.$\\
\end{cor}

In particular, in the i.i.d. standard case when $\Theta=\diag(1,\ldots,1),$  we have
$$L_{max}(x)=\begin{cases}

\frac{1}{2} (x-1)-\frac{1}{2}\log x&\textrm{if }x\ge 1,\\
\frac{r}{2} (x-1)-\frac{r}{2}\log x&\textrm{if } x\in (0,1),
\end{cases}$$
and this allows to retrieve
\cite[Corollary 2.1]{Fey}
(note that a direct proof based on the formula for the joint law of
the eigenvalues is then also available). This is in agreement with the fact that as $r$ goes to infinity, we expect
the deviations below one to be impossible
in this scale.
\\

In the general case,
we have

\begin{cor}\label{lmax}
Under the assumptions of Proposition \ref{ldpwis}, the law of the
largest  eigenvalue
of $W_n$ satisfies a large deviation principle 
in the scale $n$ with a rate function $L_{max}(\alpha)$ which satisfies, for any $\alpha \in \R,$
\begin{eqnarray*}
L_{max}(\alpha)&=& \inf\{L(\alpha_1,\ldots,\alpha_r):\max\alpha_i=\alpha\}\\
&
\ge& \inf_{\|x\|_2=1} \sup_{t\in\mathbb R}
\{t \alpha -\log E[e^{t|\langle G, \Theta^{\frac{1}{2}} x\rangle|^2}]\}
=: I_{r,\Theta}(\alpha)\\
\end{eqnarray*}
\end{cor}

From there,  one can easily  improve the upper
bound on the probability of deviations of the
largest eigenvalue of  \cite[Theorem 2.1]{Fey}~:

\begin{cor}\label{fey2.1}Assume that $G$ satisfies Assumption \ref{hyponG}
and  that the $g_i$'s are i.i.d. with mean 0 and variance 1. Let $\Theta
=\diag(\theta_1,\theta_2,\ldots, \theta_r)$ be 
a diagonal matrix with positive entries, with $\theta_1 \ge \theta_2 \ge \ldots \ge \theta_r.$
Then we have that, for $\alpha \ge \theta_1,$ 
$$ \lim_{n \ra \infty} \frac{1}{n} \log \Pro(\lambda_{\max} \ge \alpha) = -I_{r,\Theta}(\alpha).$$
\end{cor}

Note that when  $\alpha \ge \theta_1,$
 $I_{r,\Theta}(\alpha)= \inf_{[\alpha, \infty)} L_{max}$
and in particular  $I_{r,\Theta}$ is not necessarily lower semicontinuous.  
We refer to \cite{Fey} for more properties of $I_{r,\Theta}$, related results and conjectures.

\noindent
{\bf Proof of Proposition \ref{ldpwis}.}
In the case where $X_n=0$, we can apply Theorem \ref{main}
with $P_{\Theta,r}(K(z),C)=\det(z-\Theta^{\frac{1}{2}}C\Theta^{\frac{1}{2}})$ 
and ${\bf I}(C)=J(C)$. Hence, for $\alpha_1\ge \alpha_2\ge\cdots\ge\alpha_r>0$,
$L(\alpha)$ is the infimum of $J$ over the nonnegative  Hermitian matrices
$C$ such that $\Theta^{\frac{1}{2}}C\Theta^{\frac{1}{2}}$ has spectrum
$(\alpha_1,\ldots,\alpha_r)$. \qed

\noindent
{\bf Proof of Corollary \ref{ldpwisg}.}
In this case, $\log \E[e^{\langle G,Y G\rangle}]$ equals $
\log \det [ (R^{-1}-2Y)^{-\frac{1}{2}})R^{-\frac{1}{2}}]$ if $R^{-1}-2Y>0$, and is infinite otherwise.
A classical saddle point analysis shows that
the supremum in $J$ is taken at 
$$C_{ij}=\frac{ E[g_i g_j e^{\langle G,YG\rangle}]}{ E[e^{\langle G,Y G\rangle}]}= ((R^{-1}-2Y)^{-1})_{ij}$$
which yields 
$$J(C)=\frac{1}{2}\Tr(CR^{-1}-I)-\frac{1}{2} \log \det (CR^{-1}).$$
We finally take the infimum over $C$ so that $\Theta^{\frac{1}{2}}C\Theta^{\frac{1}{2}}=\sum_{i=1}^r \alpha_i e_i e_i^*$
for some orthonormal basis (ONB)
$(e_i)_{1\le i\le r}$. This gives
\begin{eqnarray*}
L(\alpha)&=&
\inf_{(e_i)ONB}
\left\{\frac{1}{2} \sum_{i=1}^r \alpha_i \langle
e_i,\Theta^{-1/2}R^{-1}\Theta^{-1/2}e_i\rangle \right\} -\frac{1}{2}\sum_{i=1}^r \log \alpha_i
-\frac{r}{2} -\frac{1}{2}\sum\log r_i(\Theta)\\
&=&\frac{1}{2} \sum_{i=1}^r \alpha_i r_i(\Theta) -\frac{1}{2}\sum_{i=1}^r \log (\alpha_i r_i(\Theta))
-\frac{r}{2} .\\
\end{eqnarray*}
\qed

\noindent
{\bf Proof of Corollary \ref{lmax}.}
We only need to take, in the definition
of $J(C)$,  $Y=t vv^*$ if $C$ has eigenvector $v$
for its largest eigenvalue to get a lower bound on $J(C)$,
and thus on $L$.\qed

\noindent
{\bf Proof of Corollary \ref{fey2.1}.} The inequality in Corollary 
\ref{lmax} gives the upper bound
and the lower bound is obtained by the same proof as in \cite{Fey}, that is by noticing that
$$\Pro(\lambda_{\max} \ge \alpha)=\Pro\left(\sup_{\|x\|_2=1}\langle x, W_n x\rangle
\ge\alpha\right)\ge  \sup_{\|x\|_2=1}\Pro(\langle x, W_n x\rangle
\ge\alpha)$$
and that for fixed $x,$ $\langle x, W_n x\rangle=n^{-1}\sum_{j=1}^n
(\langle x, \Theta^{\frac 1 2} G^j_n\rangle)^2 $
is a  sum of i.i.d. random variables so that Cramer's theorem apply.
By arguments as in \cite{Fey}, one can also check that $I_{r,\Theta}$ is increasing on $[\theta_1, \infty),$
which concludes the proof.\qed

\section{Large deviations for $H^{n}$ in the presence of outliers}\label{section_main_outliers_1_30.6.10}\mbox{}

We now go to the proof of the LDP in the presence of outliers, that will be stated in details in Theorem \ref{mainout}.
The proof follows the same lines as in the case without outliers and starts therefore with the study of the deviations of $H^n.$\\

Let  
$\mathcal K^o := \bigcup_{i=1}^{p_0} [a_i,b_i]$ a compact subset of $(b, \infty)\setminus \{\ell_1^+, \ldots, \ell_{p^+}^+\}.$
We equip again $\mathcal C(\K^o, \Hr)\times \Hr$ with the uniform topology.
Hereafter, we denote by  $\ell_i = \ell_i^+$ for $1\le i \le p^+$ and $\ell_i = \ell^-_{p^++p^--i+1}$ for $p^++1 \le i \le p^++p^-.$\\
We recall that $K^{n}(z)$ and $ C^{n}$ were defined in \eqref{defK} and \eqref{defC} respectively.

\begin{Th}\label{theo1with}
We assume that Assumptions \ref{spacing}, \ref{hyponG}, \ref{with} and \ref{Gwith} hold.
\begin{enumerate}
\item The law of 
 $\left((K^{n}(z))_{z \in \K^o}, C^{n}\right)$,
viewed as an element of
the space $\mathcal C(\K^o, \Hr)\times \Hr$ 
endowed with the uniform topology,
satisfies a 
large deviation
principle in the scale $n$ with rate function ${\bf I}^o$. For
 $K \in \mathcal C(\K^o, \Hr)$ and $C \in  \Hr$, ${\bf I}^o(K(\cdot),C)$
is infinite if $z\rightarrow K(z)$
is not uniformly Lipschitz on $\K^o$. Otherwise,
it is given by   
$$  {\bf I}^o(K(\cdot),C)= \inf\left\{ \Gamma^*(K_0(\cdot),C_0) + \sum_{i=1}^{p^++p^-} I^{(Z)}(L_i)\right\},$$
where the infimum is taken over the families  
 $K_0(\cdot)\in \mathcal C(\K^o, \Hr)$, $C_0,L_1, \ldots, L_{p^++p^-}\in \Hr$  satisfying  the condition 
\be \label{compat}
K_0(\cdot) +   \sum_{i=1}^{p^++p^-} \frac{1}{\cdot - \ell_i} L_i = K(\cdot) \textrm{ and } C_0 +   \sum_{i=1}^{p^++p^-}L_i = C
\ee
and with
\begin{multline*}
\Gamma^*(K(\cdot),C) = \sup_{P,X,Y}\left\{ \Tr\left(\int K^\prime(z)
P(z)\ud z + \sum_{i=1}^{p_0} K(b_i) X_i +C Y\right)\right. \\ \left.
-\int\Lambda\left(-\int \frac{1}{(z-x)^2}P(z)
\ud z + \sum_{i=1}^{p_0} \frac{1}{b_i-x}X_i + Y\right)\ud\mu(x)\right\},
\end{multline*}
the supremum being taken over piecewise constant $P$ with values 
 in $\Hr,$ $X= (X_1, \ldots, X_{p_0}) \in (\Hr)^{p_0}$ and $Y \in \Hr.$
\item 
The law of $(H^{n}(z))_{z \in \K^o}$ on $\mathcal C(\K^o, \mathbb R)$
equipped with the uniform topology,
satisfies a large deviation
principle in the scale $n$ with rate function
given, for a  function $f\in \mathcal C(\K^o,\mathbb R)$,
by 
$$J^o_{\K^o}(f)=\inf\{ {\bf I}^o(K(\cdot), C)\ste (K(\cdot), C)\in \mathcal C(\K^o, \Hr)\times \Hr,
P_{\Theta,r}(K(z), C)=f(z)\,\,
\forall z\in \K^o\}.
$$
\end{enumerate}

\end{Th} 

Note that the function $\Gamma^*$ is well defined because if $K$ is uniformly Lipschitz on $\K^o$, then 
so is any $K_0$ satisfying the compatibility condition \eqref{compat}, so that $K_0^\prime$ almost surely exists.\\

Under the second assertion of Assumption \ref{Gwith}, we have
the  following straightforward application
of the contraction principle.
\begin{lem}\label{LDPZ}
Let  $Z_1$ be the $\Hr$-valued random variable 
such that for $1 \le i  \le j \le r,$
$(Z_1)_{ij}=\overline{g_i(1)} g_j(1)$. Under
Assumption \ref{Gwith},
 $\frac{Z_1}{n}$   also satisfies a large deviation principle in the scale $n$
with a good rate function  $I^{(Z)}(M)=\inf\{ I(v): \overline{v_i}v_j=M_{ij}, 1\le i,j\le r\}$.
\end{lem}

The proof of Theorem \ref{theo1with}  follows the same lines 
as that of  Theorem \ref{theo1}, except that the  
 LDP for finite dimensional marginals for our process
is described by Theorem 3.2 of \cite{mylene-sandrine-jamal} instead of Theorem 2.2 of \cite{jamal-2002}. It is based
on the large deviations for 
$K^{n}$ and $C^{n}$ that can be, up to a re-indexation, shown to be exponentially equivalent to
$$K^{n}(z)_{ij}=\frac{1}{n}
\sum_{k=p^++p^-+1}^{n}\frac{1}{z-\lambda^{n}_k} 
\overline{g_i(k)} g_j(k)+ \sum_{k=1 }^{p^++p^-}\frac{1}{z-\ell_k} \frac{\overline{g_i(k)}g_j(k)}{n}
$$
which satisfy a LDP
by independence of the $g_i(k)$, and large deviations
of each parts by Proposition \ref{PGDmarginals} and Lemma \ref{LDPZ}.
The corresponding rate function will be denoted
by $(I_M^{z_1, \ldots,z_M})^o.$ To define this new rate function, we first extend in an obvious way the definition of $I_M^{z_1, \ldots,z_M}$ for $z_i$'s
 in $\K^o.$
Then one can define,
for $K_1, \ldots, K_M, C \in \Hr,$ and $z_1, \ldots, z_M \in \K^o,$
$$ (I_M^{z_1, \ldots, z_M})^o(K_1, \ldots, K_M, C )= \inf\left\{I_M^{z_1, \ldots, z_M}(K_{0,1}, \ldots, K_{0,M},C) 
+ \sum_{i=1}^{p^++p^-} I^{(Z)}(L_i) \right\},$$
where the infimum is taken over families 
$$(C,K_{0,1}, \ldots,K_{0,M}, L_1, \ldots,L_{p^++p^-})\in (\Hr)^{1+M+p^++p^-}$$  under the condition that for all $1\leq j \le M,$ 
$$ K_{0,j} +  \sum_{i=1}^{p^++p^-} \frac{1}{z_j-\ell_i} L_i = K_j \quad\textrm{ and } \quad C_0 +  \sum_{i=1}^{p^++p^-} L_i = C.$$
 By 
  Dawson-G\"artner Theorem, we deduce that $\left((K^{n}(z))_{z \in \K^o}, C^{n}\right)$ satisfies a LDP for the topology
of pointwise convergence with good rate function
$${\bf J}^o(K,C) = \sup_M \sup_{z_1< \ldots < z_M, z_j\in \K^o} (I_M^{z_1, \ldots,z_M})^o(K(z_1), \ldots, K(z_M), C).$$
Since exponential tightness is clear, this LDP can be reinforced into the uniform topology.\\
We then have to check that ${\bf I}^o={\bf J}^o.$\\

From the definition of ${\bf I}^o,$ the first thing to check is that on the event $\{{\bf J}^o(K(\cdot),C) < \infty\},$
 $K$ is  Lipschitz continuous on $\K^o$. The proof is similar
to that of Lemma \ref{proprf} as, once the $L_i$ are given, $K$ is Lipschitz
on $\K^o$ as soon as $K_0$ is.

We now suppose that $K$ is Lipschitz continuous on $\K^o$ and we want to identify the two rate functions.
By mimicking\footnote{We just have to be careful in the rewriting to put one border term for each interval involved in $\K^o.$} the proof at the end of Section \ref{subsec:proofhn}, one can easily show that for
$K$ is Lipschitz continuous on $\K^o,$ 
\be \label{Igamma}
 \sup_M \sup_{z_1, \ldots, z_M}  I_M^{z_1,\ldots, z_M}(K(z_1), \ldots, K(z_M)) = \Gamma^*(K,C).
\ee
Now, in order to achieve this identification,
 we have to check that we can switch the supremum over $M$
and the $z_i$'s and the infimum over the admissible simultaneous decompositions of $K$ and $C.$ 
It is clear that, 
$${\bf J}^o(K,C)\le \Gamma^*(K_0(\cdot),C_0) +\sum_{i=1}^{p^++p^-}I^{(Z)}(L_i)$$
for any admissible choice of $L_i$, and therefore ${\bf J}^o\le {\bf I}^o$
after optimisation.  
We now need the converse inequality. By definition of ${\bf J}^o,$ if it is finite,
then for any positive integer $p $, there exists $M(p)$ and $z_1, \ldots, z_{M(p)}$ such that
$$ {\bf J}^o(K,C)\ge (I_{M(p)}^{z_1,\ldots,z_{M(p)}})^o(K(z_1),\ldots, K(z_{M(p)}),C)-\frac{1}{p}.$$
Now for each $z_1, \ldots, z_{M(p)}$ we choose an admissible decomposition (according to \eqref{compat}) of $K$ so that
$$ {\bf J}^o(K,C)\ge I_{M(p)}^{z_1,\ldots,z_{M(p)}}(K_0^{M(p)}(z_1),\ldots, K_0^{M(p)}(z_{M(p)}),C)+\sum_{i=1}^{p^++p^-} I^{(Z)}(L^{M(p)}_i)-\frac{2}{p}.$$

Moreover, for each $M$
and choices of $z_1<\cdots<z_M$, 
$$I_M^{z_1,\ldots,z_M}(K(z_1),\ldots, K(z_M))=\Gamma^*(K^{z_1,\ldots,z_M}_M,C)$$
with $K^{z_1,\ldots,z_M}_M(z)=\sum_{i=1}^M 1_{[z_i,z_{i+1}]} K(z)$.

 By definition, since $I^{(Z)}$ and $\Gamma^*$
are good rate functions and as for all $i,$ $I^{(Z)}(L_i^{M(p)})$ and $\Gamma^*(K_0^{M(p)}(z_1),\ldots, K_0^{M(p)}(z_{M(p)}),C)$
are uniformly bounded, it  implies that the arguments are tight and we can take a converging subsequence. 
Let $K_0$ and $L_i$ be limits along a subsequence, we get
$${\bf J}^o(K,C)\ge  \Gamma^*(K_0(\cdot),C)+\sum_{i=1}^{p^++p^-} I^{(Z)}(L^i)
-\frac{1}{p}$$
which insures that ${\bf J}^o(K,C)\ge {\bf I}^o(K,C)$. 
This completes the proof of Theorem \ref{theo1with}.\\

\section{Large deviations principle 
for the largest eigenvalues 
in the case with outliers}\label{section_main_outliers_30.6.10}\mbox{}

We now state the main theorem of this section,
 namely an analogue of Theorem \ref{main}.
For any $\e, \rho $ small enough, we define the compact sets
$$\K^o_{\e,\rho}:= [b+\e,\e^{-1}] \setminus \bigcup_{i=1}^{p^+} (\ell_i^+-\rho, \ell_i^++\rho)$$
and $\K^o_{\e}:= [b+\e,\e^{-1}].$
We also define the set $\{\ell\}:= \{\ell_1^+, \ldots, \ell_{p^+}^+,b\}, $ and
for $z \notin \{\ell\},$  $R(z):= \prod_{i=1}^{p^+} \frac{1}{z-\ell_i^+}.$
We recall that $s$ is the sign of the product $\prod_{i=1}^r \theta_i.$ \\

For any $\e,\rho, \gamma > 0,$ and $\alpha \in \R^p_\downarrow (b+\e),$ we put 
\begin{multline*}
S_{ \alpha, \gamma}^{\e,\rho,o}:=\left\{f\in \mathcal C(\K^o_{\e,\rho},\R):  \exists
 g \in \mathcal C(\K^o_{\e,\rho},\R) \textrm{ with }  \gamma\le g\le \frac{1}{\gamma}, \omega(g) \le \frac{1}{\gamma} \right.\\
 \left. \textrm{ and }
f(z)=s.R(z).g(z) \prod_{i=1}^p (z-\alpha_i)\right\}
\end{multline*}
We also denote by
\begin{multline*}
C_{k,\gamma}^{\e,\rho,o}:= \left\{f\in \mathcal C(\K^o_{\e,\rho},\R): \exists p\textrm{ polynomial of
degree $m+p^+-k$ with } m+p^+-k \textrm{ roots in } \K^o_\e \right.\\
 \left.\textrm{ and dominant coefficient } 1,
 g \in \mathcal C(\K^o_{\e,\rho},\R) \textrm{ with }    \gamma\le g\le \frac{1}{\gamma}, \omega(g) \le \frac{1}{\gamma} \right.\\ \left. \textrm{ and }
f(z)=s.g(z).R(z).p(z)\right\}
\end{multline*}
and 
$$C_{\gamma}^{\e,\rho,o} = \bigcup_{0 \le k \le m+p^+ }C_{k,\gamma}^{\e,\rho,o}.$$

Then the main statement of this section is the following.

\begin{Th}\label{mainout}
Under Assumptions \ref{spacing}, \ref{hyponG}, \ref{with} and \ref{Gwith},
the law of the $m+p^+$ largest eigenvalues $(\wtl_1^{n},\ldots, \wtl_{m+p^+}^{n})$ of $\wtX$  satisfies
a large deviation principle in $\R^{m+p^+}$ with good rate function $ L^o.$ 
For $\mathbf \alpha = (\alpha_1, \ldots, \alpha_{m+p^+})\in \R^{m+p^+},$ 
we take $\alpha_{m+p^++1} =b$ and  $L^o$ is defined as follows :
\[
L^o( \alpha ) = \left\{
\begin{array}{ll}
\lim_{\e \downarrow 0}\lim_{\rho \downarrow 0}\inf_{\cup_{\gamma >0} S_{(\alpha_1, \ldots, \alpha_{m+p^+-k}), \gamma}^{\e,\rho,o}} J^o_{\K_{\e,\rho}^o} &
\textrm{if } \alpha \in \R^{m+p^+}_\downarrow(b), \alpha_{m+p^+-k+1} = b, \\
& \textrm{ } \alpha_{m+p^+-k} >b
  \textrm{ for a }
k\in \{0, \ldots, m\}, \\ \\
\infty & \textrm{otherwise.}
\end{array}
\right.
\]
\end{Th}

Even though the rate function $L^o$ is not very explicit, we show below that it must be infinite if Horn's inequalities are violated.

\begin{rmq}
Recall that the eigenvalues $(\widetilde\lambda^{n}_i)_{1\le i\le n}$
 of the sum of two Hermitian
matrices  with eigenvalues $(\lambda^{n}_i)_{1\le i\le n}$ 
and $\theta:=(\theta_1, \ldots, \theta_{ r},0,\ldots,0)$
 satisfy   Horn's inequalities and are characterised by the
fact that they satisfy such  inequalities (see \cite{TK} for details).
Assume that $\wtl:=(\wtl_1,\ldots,\wtl_{m+p^+})$ is at distance
of  the bulk and of 
the outliers which is bounded below.  We claim that 
 the rate function $L^o(\wtl)$ is  infinite if $(\wtl, \ell, \theta) $
do not satisfy the Horn inequalities. Indeed, if $L^o(\wtl)$ is finite,
$(\wtl_1,\ldots,\wtl_{m+p^+})$ are  zeroes of a function $f$
which can be written
$$f(z)=P_{\Theta,r}(K(z),C).$$
with
${\bf I}^o(K(\cdot),C)$  finite. It implies that
 there exists  sequences $\lambda^{n}\in\mathbb R^n, g_j(\cdot)\in \mathbb C^n$
so that  $\lambda^{n}$ satisfies Assumptions  \ref{spacing},  \ref{with} and \ref{Gwith}
and 
 $$K^{n}(z)=\frac{1}{n}\sum_{i=1}^n \frac{\ovl{ g_i(k)}
g_j(k)}{z-\lambda^{n}_k}
,\quad C^{n}=\frac{1}{n}\sum_{i=1}^n \ovl{ g_i(k)}
g_j(k)$$
converge to $K(z)$ (uniformly away from the bulk and the
outliers) and $C$ respectively.
By definition, there exists a constant $c$ such that
$$P_{\Theta,r}(K^{n}(z),C^{n})\prod_{i=1}^n
 (z-\lambda^{n}_i)=c\det\left(z-\diag(\lambda^{n})-\sum_{i=1}^r
 \theta_iu_i u_i^*\right)$$
with $u_i=g_i$ in the i.i.d.  perturbation model and $u_i$ the Gram-Schmidt orthonormalization 
of the vectors $g_i$ in the orthonormalized  perturbation model. Hence,
 the function  $f_n(z)=P_{\Theta,r}(K^{n}(z),C^{n})$
vanishes at the eigenvalues $(\widetilde\lambda^{n})$
 of the sum of the two
Hermitian matrices $\diag(\lambda^{n})$ and $\sum_{i=1}^r
 \theta_iu_i u_i^*$ (note that we can assume without
loss of generality that its zeroes are different from $\lambda^{n}$
by Lemma \ref{vpcommon}). Therefore, $(\lambda^{n},\wtl^{n},
\theta)$ satisfy Horn's inequalities by \cite{TK}.
Since the $(\wtl^{n})$ are bounded, they are
relatively compact and we see that the limit points $(\widetilde\lambda_1,
\ldots,\wtl_{m+p^+})$ of $(\wtl_1^{n},\ldots,\wtl_{m+p^+}^{n})$ which stay away from the bulk and the outliers
are the zeroes of $f$. By passing to the limit in Horn's inequalities,
we thus deduce that if the vector  $(\wtl_1,\ldots,\wtl_{m+p^+})$ has finite $L^o$-entropy,
and is 
 away from the bulk and the outliers,  $(\wtl, \ell, \theta) $ satisfies Horn's inequalities. 
It would be interesting to have a direct proof
of this fact.
 
\end{rmq}

\subsection{Proof of Theorem \ref{mainout}}\mbox{}\\

We now prove Theorem \ref{mainout}, following roughly
 the same lines as for Theorem \ref{main}.

As in  the proof of Theorem \ref{main}, the crucial point is to use Proposition \ref{zeroesscheme}.
In the sticking case, if $z \in \K_\e,$ for $n$ large enough,
the condition that $z$ should not belong to the set of eigenvalues of $X_n$ was very easy to check.
Here, we need to make sure that the eigenvalues
are not exactly equal to the outliers
to use our strategy. We  show the following

\begin{lem} \label{lemdistinct}
Assume that the eigenvalues $\la_1^{n}, \ldots, \la_n^{n}$ of 
$X_n$ are pairwise distinct and that  Assumptions \ref{spacing} and \ref{with} hold, then $X_n$ and $\wtX$ have no eigenvalue in common
for almost all $G$.
\end{lem}

The proof of this lemma is postponed to Appendix \ref{SectionAppendixdistinct}.
We shall therefore give the proof of the Theorem when the eigenvalues
of $X_n$ are distinct. This is however sufficient to get 
the LDP without this hypothesis 
due to
 the following Lemma.

\begin{lem}\label{separating}
Let $X_n$  satisfy Assumptions \ref{spacing} and \ref{with}. 
Then,   there exists a sequence  $\bar X_n$ of matrices with pairwise distinct eigenvalues   satisfying Assumptions  \ref{spacing} and \ref{with}  such that,
if we define $\widetilde {\bar X_n}$ be the perturbation of $\bar X_n$
by the i.i.d. or the orthonormalized vectors constructed on the
law $\mu_n=\mu*\gamma_n$ of $G+\e(n)A$ with $A$ $r$ independent standard nornal variables
 and $\varepsilon(n)$ going to zero with $n$ fast enough, 
then, with  $(\widetilde{\bar{\lambda}^{n}_{i}})_{i\le m}$
 the extreme eigenvalues of $\widetilde {\bar X_n},$
$$\limsup_{n\ra\infty}\frac{1}{n}\log P\left(
\max_{1\le i\le m}| \widetilde{\lambda^{n}_i}-\widetilde{\bar\lambda^{n}_i}|\ge
\frac 1 n\right)=-\infty.$$
\end{lem}

\begin{pr}
We take 
$\bar X_n$ to be the matrix with the same eigenvectors as $X_n$ 
and the same eigenvalues except for those which are sticked
together which we separate by an arbitrary small weight $w_n\le 1/n$, much
smaller than the minimal distance
between two distinct eigenvalues of $X_n$, 
so that the eigenvalues of $\bar X_n$
are distinct and the operator norm of $X_n-\bar X_n$ is bounded above 
by $w_n$. It is straightforward to verify 
 Assumptions \ref{spacing} and \ref{with} for $\bar{X}_n$.
Now, if we add the same perturbation to $X_n$ and $\bar{X}_n$ respectively,
their eigenvalues will differ at most by $w_n$ almost surely.  Then
adding a Gaussian vector of variance $\e(n)^2$ to $G$
will not change the eigenvalues by more than $\sqrt{\e(n)}$
with probability greater than $1-e^{-\e(n)^{-1} n}$ as the empirical 
covariance matrix of this additional term is bounded by $C \sqrt{\e(n)}$ with
such a probability. We conclude by choosing $\e(n)$ such that $\sqrt{\e(n)}<1/n.$
\end{pr}

Lemma \ref{separating} means in particular the random variables $(\widetilde{\bar{\lambda}^{n}_{i}})_{i\le m}$
and  $({ \widetilde\lambda^{n}}_i)_{i\le m}$ are exponentially equivalent and
 \cite[Theorem 4.2.13]{DZ} asserts that a large deviations principle for
the extreme eigenvalues 
$(\widetilde{\bar{{\lambda}^{n}_{i}}})_{i\le m}$ 
 of $\widetilde{\bar {X_n}}$ entails the large deviations principle for the law of  $({ \widetilde{\lambda^{n}_i}})_{i\le m}$
with the same rate function. Therefore, the proof of Theorem \ref{with} can be done for the eigenvalues of $\widetilde {\bar X_n},$
the main advantage being that, from Lemma \ref{lemdistinct} above, we get that $\bar X_n$ and $\widetilde {\bar X_n}$
have almost surely no eigenvalue in common and we can proceed
as in the  case without outliers.\\

From now on, we assume that $X_n$ satisfies Assumptions \ref{spacing} and \ref{with} and has pairwise distinct eigenvalues
and that $G$ satisfies Assumptions \ref{hyponG} and \ref{Gwith} and that its law is absolutely continuous with respect to Lebesgue measure.\\

We first  focus our attention 
to the function $H^{n}$ restricted to $\K_{\e,\rho}^o$ and show the counterpart of Lemma \ref{lemonH}, that is
\begin{lem} \label{hh}
Let $\e,\rho $ be fixed. There exists  a positive integer $n_0(\e, \rho)$ and  $L(\e) >0$ such that for any $n \ge n_0(\e,\rho),$
for any $z \in \K_{\e,\rho}^o,$
 \begin{equation} 
H^{n}(z)=\begin{cases}s\prod_{i=1}^r \|q_i^{n}W_i^{n}\|_2^2 \,\,
g_n(z)R(z)\,\,
\prod_{i=1}^m (z-\wtl_i^{n})&\textrm{in the orth.  perturb. model,}\\
&\\
s \,\,
g_n(z)R(z)\,\,
\prod_{i=1}^m (z-\wtl_i^{n})&\textrm{in the i.i.d.  perturb. model,}\end{cases}
\end{equation}
with $ L(\e)\le g_n \le \frac{1}{L(\e)}$  and $\omega(g_n) \le \frac{1}{L(\e)}.$

In particular, for any $\e >0$ and $\rho >0$ small enough,
$$
\limsup_{\gamma\downarrow 0}
\limsup_{n\ra\infty}\frac{1}{n}\log \Pro\left( \left(\prod_{i=1}^m (z-\wtl_i^{n})^{-1}
H^{n}(z)\right)_{z \in 
\K_{\e,\rho}^o} \in (C_{\gamma}^{\e,\rho,o})^c\right)= -\infty.
$$
\end{lem}

\begin{pr} 
In this case,
$$ g_n(z):= \prod_{i=1}^r |\theta_i|^{-1}\prod_{i=1}^{p^+} \left(1+ \frac{\lambda_i^{n}- \ell_i^+}{z-\lambda_i^{n}}\right)
\prod_{i=p^++1}^{m+p^+} \frac{1}{z-\lambda_i^{n}}\prod_{i=m+p^++1}^n 
 \left(1+
\frac{\lambda_i^{n} -\wtl_i^{n}}{ z-\lambda_i^{n}}
\right).$$

The proof is exactly the same as in the sticking case once we have noticed that,
from Assumption \ref{with}, there exists $n_0(\e,\rho)$ such that for $n \ge n_0(\e,\rho),$  
$\prod_{i=1}^{p^+} \left(1- \frac{\lambda_i^{n}- \ell_i^+}{\lambda_i^{n}- z}\right) \ge \frac{1}{2^{p^+}},$
so that
$$ g_n(z) \ge \prod_{i=1}^r |\theta_i|^{-1}\left(\frac{1}{2}\right)^{p^+}\left(\frac{\e}{2}\right)^{m} e^{-2m\left(\frac{4}{\e^2}\right)^2}.$$
Note that we could similarly show that for $n\ge n_0(\e,\rho),$ 
\be \label{upperwith}
 g_n(z) \le L(\e):=\left(\frac{3}{2}\right)^{p^+}\left(\frac{2}{\e}\right)^{m} e^{\frac{8m}{\e^2}}.
\ee
The uniform equicontinuity is also shown very similarly.
\end{pr}

As in the sticking case,
we have  the analogue of Lemma \ref{near}, with $L^o$ instead of $L.$
To state more precisely the lemma,  we introduce the following notation:
we denote by 
$ G_k(\alpha, \delta, \e, \rho) $ the set of $n$tuples $(\wtl_1^{n}\ge  \cdots\ge \wtl_n^{n})$  \st for all $i\le m+p^+-k$, 
$$|\wtl^{n}_i-\alpha_{i}|\le\delta\textrm{ if $ \alpha_i \notin \{\ell\}\qquad$ and $\qquad$}|\wtl^{n}_i-\alpha_{i}|\le\rho\textrm{ if $\alpha_i \in \{\ell\}$}$$ and for all ${m+p^+-k+1\le i\le m+p^+}$, 
$$\wtl^{n}_i\le b+\e.  $$
Because of Lemma \ref{hh}, $H_n$ belong 
to the set of functions $f(z)=h(z)\prod_{i=1}^m (z-\alpha_i)
R(z)$  with a bounded positive constant $h$ on $\K^o_{\e,\rho}$
with values in $[\gamma,\gamma^{-1}]$
 with overwhelming probability.
But on this set also the zeroes $\alpha_i$ are continuous function
of the functions $f$ and therefore we can proceed
exactly as in the case without outliers.

\begin{lem} \label{nearwith}
Let $\alpha \in \R^m_\downarrow$ and $k$ between $0$ and $ m$ such that $\alpha_{m+p^+-k+1} = \ldots = \alpha_{m+p^++1}=b$ and $\alpha_{m+p^+-k}>b.$
We have
\begin{multline*}
\lim_{\e \downarrow 0}\lim_{\rho \downarrow 0}\limsup_{\delta \downarrow 0}
\limsup_{n\ra\infty}\frac{1}{n}\log \Pro\left( (\wtl_1^{n}, \ldots, \wtl_n^{n}) \in G_k(\alpha, \delta, \e, \rho)\right)
\\
=\lim_{\e \downarrow 0}\lim_{\rho \downarrow 0}\liminf_{\delta \downarrow 0}
\liminf_{n\ra\infty}\frac{1}{n}\log \Pro\left( (\wtl_1^{n}, \ldots, \wtl_n^{n}) \in G_k(\alpha, \delta, \e, \rho)\right)
=-L^o(\alpha),
\end{multline*}
with the obvious convention that $\bigcap_{m+p^+-k+1\le i\le m+p^+}\{\wtl^{n}_i\le b+\e \}= \Omega$ if $k=0.$
\end{lem}

The proof is similar to the case without outliers.

\section{Application to $X_n$ random,
following some classical matrix distribution}\label{Proof_classical_300610}
This section is devoted to the proofs of the results stated in Section \ref{secmm}.
\subsection{Proof of Theorem \ref{tholdx}} \mbox{}\\
Theorem \ref{tholdx} is a slight extension of \cite[Theorem 2.6.6]{alice-greg-ofer} and the proof will therefore follow the same lines.
We introduce the notations $\phi(\mu,x)=-V(x)+\beta \int \log|x-y|d\mu(y)$ (for $x$ greater or equal the right edge of the support of $\mu$)
and $\hat\mu_n=\frac{1}{n-p}\sum_{i=p+1}^n \delta_{\lambda_i^{n}}$. Then
\begin{multline*}
 \Pro_{V,\beta}^n(\ud\lambda_1,\ldots,\ud\lambda_n)=\\
\frac{
Z_{nV/(n-p),\beta}^{n-p}
}{Z_{V,\beta}^n}
  e^{n\sum_{i=1}^p \phi(\hat\mu^n,\lambda_i)+ \beta
\sum_{1\le i<j\le p}\log |\lambda_i-\lambda_j|} \ud \Pro^{n-p}_{n V/(n-p), \beta}
(\lambda_{p+1},\ldots, \la_n)\ud\lambda_1\cdots\ud\lambda_p.
\end{multline*}
By \cite[Lemma 2.6.7]{alice-greg-ofer}, if parts $i)$ and $ii)$ of Assumption \ref{assX} hold, the law $ \Pro_{V,\beta}^n$
is exponentially tight so that it is enough to estimate the probability of a small ball around ${\bf x}=(x_1\ge x_2\ge\cdots\ge x_p)$
(with $x_p \ge b_V$),
namely events of the form $B({\bf x},\delta) : = \{\max_{1\le i\le p}
|\lambda_i-x_i|\le \delta, \max |\lambda_i| \le M\}.$\\

As in  \cite{bdg01}, a crucial observation is
 the fact that $\hat\mu_n$ converges
to $\mu_V$ much faster than exponentially under $\Pro^{n-p}_{n V/(n-p), \beta}$ (its LDP is indeed in the scale $n^2$).
We can  therefore replace $ \phi(\hat\mu^n,\lambda_i)$
by $\phi(\mu_V, x_i)$,  whereas 
the ratio of partition functions converges by hypothesis.\\

To be more precise, let us first sketch the proof of  the upper bound. 
Note that there exists a constant $\Phi_M$ such that on $B({\bf x},\delta),$ $\phi$
is bounded above by $\Phi_M$  so that
\begin{eqnarray*}
 \Pro_{V,\beta}^n(B({\bf x},\delta))&\le& \frac{Z_{nV/(n-p),\beta}^{n-p}}{Z_{V,\beta}^n} e^{\beta p(p-1)/2\log (x_1-b_V)}
(e^{np \Phi_M}\Pro^{n-p}_{n V/(n-p), \beta}(\hat\mu^n\in B_\varepsilon(\mu_V)^c)
\\
&&+  (2M)^p e^{n\sum_{i=1}^p \max_{|y-x_i|\le\delta}
\max_{\mu\in B_\varepsilon(\mu_V)}\phi(\mu, y)})
\end{eqnarray*}
with $ B_\varepsilon(\mu_V)$ a small ball with radius $\varepsilon$ around
$\mu_V$ for a distance compatible with the weak topology. 
As $n^{-2}\log \Pro^{n-p}_{n V/(n-p), \beta}(\hat\mu^n\in B_\varepsilon(\mu_V)^c)$
is bounded above by a negative real number for all $\delta>0$
by the LDP for the law of $\hat\mu^n$, the first term is negligible
as $n$ goes to infinity.
 Using the fact that $(\mu,x)\ra \phi(\mu,x)$ is upper continuous, we obtain  the upper bound by first letting $n$ go to infinity, then letting $\delta$ decrease to
zero and finally letting $\e$ go to zero.  Notice again that in the proof of this upper bound, 
we use part $i)$ of Assumption \ref{assX} to get the LDP for $\hat\mu_n$ and $ii)$ to control the ratio of the partition functions.\\

The lower bound is similar to the proof in \cite[p. 84]{alice-greg-ofer}, which corresponds to $p=1.$
We proceed by induction on $p$ and we can therefore assume that $p$ is the smallest integer
such that  $x_p > b_V.$
There exists $x_i^\delta, 1\le i\le p$, whose
small neighbourhood are included 
in the $\delta$ neighbourhood of
$x_i,1\le i\le p$, and  which are distinct, so that for $\varepsilon$
small enough
\begin{multline*}
 \Pro_{V,\beta}^n(\max_{1\le i\le p}
|\lambda_i-x_i|< \delta)\ge \Pro_{V,\beta}^n(\max_{1\le i\le p}
|\lambda_i-x_i^\delta|< \varepsilon, \la_i < x_p - \delta - \varepsilon, \forall i >p)\\
\ge \frac{Z_{nV/(n-p),\beta}^{n-p}}{Z_{V,\beta}^n} \exp\bigg( (n-p) \inf_{\tiny \begin{array}{cc}
|y_i-x_i^\delta|< \varepsilon\\
\mu \in B_{[-M,x_p-\delta-\varepsilon]}(\mu_V, \varepsilon) \end{array}}
 \phi(y_i, \mu)\bigg) \Pro^{n-p}_{n V/(n-p)}(\hat\mu^n\in  B_{[-M,x_p-\delta-\varepsilon]}(\mu_V, \varepsilon)),
\end{multline*}
with $B_{[-M,x_p-\delta-\varepsilon]}(\mu_V, \varepsilon))$ the set of probability measures in $B_\varepsilon(\mu_V)$
with support in $[-M,x_p-\delta-\varepsilon].$ 
When the $x_i$'s are distinct and away from $b_V,$
their logarithmic interaction is negligible; moreover, part $iii)$ of Assumption \ref{assX}
allows to claim that the last term in the lower bound above converges to one.
 We therefore get
$$\liminf_{n\ra\infty}\frac{1}{n}\log \Pro_{V,\beta}^n(\max_{1\le i\le p}
|\lambda_i-x_i|< \delta)\ge -\sum_{i=1}^p J_V(x_i^\delta)
-\alpha_{V,\beta}^p$$
Now, $J_V$ is continuous away from the support of $\mu_V$
so that we can conclude by letting $\delta$ going to zero. \\
Then to get the correct expression of the rate function, we just have to check that $\alpha_{V,\beta}^p=p\alpha_{V,\beta}^1,$
which is easy and left to the reader.
\qed

\subsection{Proof of Theorem \ref{theorandom}}\mbox{}\\
As explained in Section \ref{secmm}, we have to study $\wtX,$ when $X_n$ is diagonal
with eigenvalues having $\Pro_{V, \beta}^n$ as their joint law and the $U_i$'s obtained
by orthonormalisation procedure from $G= (g_1, \ldots, g_r)$ i.i.d. standard Gaussian.\\ 
The proof will consist in first
fixing the possible deviations of the extreme eigenvalues of
$X_n$ (hence providing outliers) and then, being given these outliers,
computing the deviations of the eigenvalues of $\wtX$.
The main point of course is that
with exponentially large probability, only
a finite number of eigenvalues of $X_n$ can deviate.

\noindent
\noindent
$\bullet$ More precisely,  we observe that, by Theorem \ref{tholdx}, 
for all $p \in \n^*$, the probability
that $\lambda_p^{n}$ is greater  than $b_V+\delta$
is less than $e^{-n p\varepsilon(\delta)}$
for $\varepsilon(\delta)=\inf_{(b_V+\delta,+\infty)}J_V.$ The only point to check is that
 $\inf_{(b_V+\delta,+\infty)}J_V > 0,$ which is a consequence of part $iii)$ of 
Assumption \ref{assX}.\\

\noindent
$\bullet$ The deviations of the eigenvalues of $X_n$ are controlled by Theorem \ref{tholdx} :  there exists 
$\e(\eta,\ell)>0$ so that  for $n$ large enough,
for $\e\le\e(\eta,\ell)$
\begin{equation}\label{ldpes}-J^p(\ell_1,\ldots,\ell_p)-\eta\le 
\frac{1}{n}\log \Pro\left( \max_{1\le i\le p} |\lambda_i^{n}-\ell_i|\le
\e,  \lambda_{p+1}^{n}\le b_V+\e\right)
\le-J^p(\ell_1,\ldots,\ell_p)+\eta.\end{equation}

\noindent
For all $
(\ell_1,\ldots,\ell_p)$ and $\eta >0,$ we
define the set \beq &V_\eta(\ell_1,\ldots,\ell_p)=\left\{(\la_1, \ldots ,\la_n)\in \R^n\ste  
  \max_{1\le i\le p} |\lambda_i-\ell_i|<
\e(\eta,\ell),  \lambda_{p+1}< b_V+\e(\eta,\ell)\right\}.&\eeq

$\bullet$ Now, knowing  the deviations of the eigenvalues of $X_n,$ one can treat them as outliers
and deal with the eigenvalues of the perturbed model. We have that,
for any $(\ell_1,\ldots,\ell_p)\in (b_V,+\infty)^p$ and
any $\eta>0,$ there exists
$\e(\eta,\ell),\delta(\eta,\ell)>0$ so that for $n$ large enough,
for $\e<\e(\eta,\ell),\;\delta<\delta(\eta,\ell)$,
\begin{equation}\label{lkj}-L^0_{\ell_1,\ldots,\ell_p}
(x_1,\ldots,x_k)-\eta
\le 
\frac{1}{n}\log \Pro\left( \max_{1\le i\le k} |\wtl_i^{n}-x_i|\le
\delta\bigg|{\max_{1\le i\le p} |\lambda_i^{n}-\ell_i|\le
\e,\atop \lambda_{p+1}^{n}\le b_V+\e}\right)\quad\end{equation}
$$\qquad\qquad\le -L^0_{\ell_1,\ldots,\ell_p}
(x_1,\ldots,x_k)+\eta$$
These inequalities are a consequence of Theorem \ref{mainout}.
Indeed, let $X_n$ be a matrix such that the event $\{\max_{1\le i\le p} |\la_i^{n}-\ell_i|\le
\epsilon\}$ holds.
Let  $X_n'$ be a real diagonal matrix with same eigenvalues of $X_n$ except its $k$
largest eigenvalues are equal to the outliers   $(\ell_1,\ldots,\ell_p).$
Then we have $\|X_n-X_n'\|_\infty\le \epsilon$, so that, with obvious 
notations,  $\|\wtX-\wtX'\|_\infty\le \epsilon$, so that the ordered eigenvalues
of $\wtX$ and $\wtX'$ differ at most by 
$\epsilon$. Thus,  up to change $\delta$ into $\delta\pm\e$,  Theorem \ref{mainout} gives \eqref{lkj}.\\

\noindent
$\bullet$ We have now all the ingredients to prove the LDP. It is clear that since the largest eigenvalues of $X_n$ are
exponentially tight, so are the eigenvalues 
of $\wtX$, and therefore it is enough
to prove a weak large deviation principle.
 We let $K(L)$ be such
that the probability that $\lambda_1^{n}$ or $\wtl_1^{n}$
is greater than $K(L)$ is smaller than $e^{-nL}$. 

\noindent
$\bullet$ To prove the upper bound we can write, for any $p \ge k,$ any $\eta >0,$ $\delta >0,$
\begin{eqnarray}
 \Pro\left( \max_{1\le i\le k} |\wtl_i^{n}-x_i|\le
\delta\right)&\le& \Pro\left( \max_{1\le i\le k} |\wtl_i^{n}-x_i|\le
\delta, \; \lambda_{p+1}^{n}\le b_V+\delta\right) +
e^{-n p\varepsilon(\delta)}
\label{bnm}
\end{eqnarray}
We fix $\eta >0.$
As $[b_V,K(L)]^p$ is compact, from   its infinite open  covering
$\cup V_\eta(\ell_1,\ldots,\ell_p),$ one can always extract a finite covering $\cup_{1\le s\le M(\eta)} V_\eta(\ell_1^s,\ldots,\ell_p^s).$ We then take $\delta=\min \delta(\eta,\ell^s)>0$.
 Thus, we get by the LDP estimate \eqref{ldpes}

\begin{eqnarray*}
\Pro\bigg( \max_{1\le i\le k} |\wtl_i^{n}-x_i|\le
2\delta\bigg)& \le &e^{-nL}+e^{-n p\varepsilon(\delta)} +\sum_{s=1}^{M(\eta)}
\Pro\left(  \max_{1\le i\le k} |\wtl_i^{n}-x_i|\le
2\delta \cap
 V(\ell_1^s,\ldots,\ell_p^s)\right)\\
&\le &e^{-nL}+e^{-n p\varepsilon(\delta)}+\sum_{s=1}^{M(\eta)} e^{-n L^0_{\ell_1^s,\ldots,\ell_p^s}
(x_1,\ldots,x_k)
-nJ^p(\ell_1^s,\ldots,\ell_p^s)+n\eta}\\
&\le& M(\eta)e^{-n \min_{1\le s\le M(\eta)} 
(L^0_{\ell_1^s,\ldots,\ell_p^s}(x_1,\ldots,x_k)+
J^p(\ell_1^s,\ldots,\ell_p^s)-\eta)}
+ e^{-nL}+e^{-n p\varepsilon(\delta)}\\
&\le&
3M(\eta)e^{-n\min\{ L, p\varepsilon(\delta),  \tilde J^k(x_1,\ldots,x_k)\}
}\\
\end{eqnarray*} 
which gives the announced bound by taking 
first the limit as $n$ goes to infinity, then $L,p$ to infinity
and finally  $\delta$ and $\eta$ to zero.

$\bullet$ The lower bound is easier as we simply write
$$\Pro\left( \max_{1\le i\le k} |\wtl_i^{n}-x_i|\le
2\delta\right)
\ge \Pro\left(  \max_{1\le i\le k} |\wtl_i^{n}-x_i|\le
2\delta\cap 
 V(\ell_1^s,\ldots,\ell_p^s)\right)$$
and use the large deviation theorems. \qed

\section{Appendix}
\label{SectionAppendix}

\subsection{Proof of a technical lemma}\mbox{}\\

With the notations of Section \ref{SectionGS}, we have the following result

\begin{lem}\label{norm}
Under Assumption \ref{hyponG},
 for any $1\le i_0 \leq r$,  we have   
$$ \lim_{\delta \downarrow 0} \limsup_{n \ra \infty} \frac{1}{n} \log \Pro\left(\|q_{i_0}^{n}W_{i_0}^{n}\|_2^2 \notin \left[ \delta,\frac{1}{\delta}\right]\right) = - \infty.$$
\end{lem}

\begin{pr}
To simplify the notations,
we shall assume that $i_0=r$.

  Recall that the $G_i^{n}$'s were constructed from a family  $(G(k)=(g_1(k), \ldots, g_r(k))_{\substack{ k\ge 1}}$ of independent copies of $G$,
 via the formula  $G_i^{n}:=(g_i(1), \ldots,g_i(n) )^T$. For $1 \le k $, we consider the random $r\times r$ Hermitian matrix 
 $$ Z_k = G(k)^*G(k)= [\ovl{g_i(k)}g_j(k)]_{1\le i,j \le r} \quad \textrm{ and } \quad  L^{n} = \frac{1}{n } \sum_{k=1}^n Z_k.$$
  
By Cram\'er's Theorem \cite{DZ}, we have that the law of $L^{n}$ satisfies a LDP
with convex good rate function 
$$ I^{(L)} (y) =\sup_{\la\in {\Hr}}\{\lan \la, y\ran-\Lambda (\la)\},$$  where 
$ \Lambda (\la) = \log \E(e^{\langle \lambda,Z_1\ran })$ is exactly the function defined in Equation \eqref{1.02.10.def_La}.

Note that since for all $n$, $L^{n}$ is almost surely a positive semi-definite matrix, by closedness of the set of such matrices,
 the domain of $I$ is contained in the set of positive semi-definite matrices. 

 Let $P_r$ be the real polynomial function on ${\Hr}$ introduced in Proposition \ref{pol}: 
we have $\|q_{r}^{n}W_{r}^{n}\|_2^2 = P_{r}(L^{n} )$.  
Therefore, if, for any $\delta >0,$ we introduce the closed set $\mathcal E_\delta := \{y \in {\Hr}\ste 
P_{r}(y) \leq \delta \}$, we have
$$ \limsup_{n \ra \infty} \frac{1}{n} \log \Pro(\|q_{r}^{n}W_{r}^{n}\|_2^2 \leq \delta)  \leq - \inf_{y \in \mathcal E_\delta } I^{(L)} (y).$$
Let us assume that 
$$ M:= \lim_{\delta \downarrow 0} \inf_{y \in \mathcal E_\delta } I^{(L)} (y)< \infty.$$

Since $I^{(L)}$ is a good rate function, there exists a compact set $K$ such that 
$\inf_{y \in K^c }I^{(L)}(y) >M$, so that for all $\delta>0$,  $\inf_{y \in \mathcal E_\delta }I^{(L)}(y)= 
\inf_{y \in \mathcal E_\delta \cap K} I^{(L)} (y).$
Moreover the infimum on $\mathcal E_\delta$  is reached : 
  let, for all $n\ge 0$,  $y_n$ be an element of $K$ such that   $I^{(L)}(y_n)= \inf_{y \in \mathcal E_{\frac{1}{n}} }I^{(L)}(y).$ \\
There exists a subsequence $\varphi(n)$ such that $y_{\varphi(n)}$ converges, as $n$ goes to infinity to some $y_0.$ 
By continuity of $P_{r}$, $P_{r}(y_0) = \lim_{n \ra \infty} P_{r}(y_{\varphi(n)}) = 0$.
It follows, by the last part of Proposition \ref{pol}, that $y_0$ is not positive definite. 
However, since $I^{(L)}$ is lower semicontinuous, we have
 $I^{(L)}(y_0) \leq M < \infty$, which implies that  $y_0$ is a positive semi-definite matrix. Let $p$ be the orthogonal projection onto $\ker y_0$. Note that $p\neq 0$ and that $\lan p,y_0\ran=\Tr (y_0p)=0$. 
 \beq I^{(L)}(y_0)&=&\sup_{\la\in {\Hr}}\{\lan\la,y_0\ran-\La(\la)\}\\
&\geq& \sup_{t>0}\{ \langle-t p, y_0 \rangle -\La(-tp)\}\\
&=&\sup_{t>0} - \La(-tp)\\
&=&+\infty \qquad\textrm{by  \eqref{La_inf_inf.1.02.10},}
\eeq 
which yields a contradiction (as  we already proved that $I^{(L)}(y_0)\leq M$).

Similarly, as $I^{(L)}$ is a good rate function,
it has compact level sets and therefore has to be large
on the set $\{y: P_r(y)\ge 1/\delta\}$. Hence,
$$\limsup_{\delta\downarrow 0}
 \limsup_{n \ra \infty} \frac{1}{n} \log \Pro(\|q_{r}^{n}W_{r}^{n}\|_2^2 \geq \delta^{-1})  =-\infty$$
which completes the proof of the lemma.
 \end{pr}
 
\subsection{On the eigenvalues of the deformed matrix} \mbox{}\\
\label{SectionAppendixdistinct}

The goal of this section is to prove Lemma \ref{lemdistinct}. In fact,
we will prove the slightly more general

\begin{lem} \label{vpcommon}
 Let $\mathbb{K}$ be either $\R$  or $\C.$ Let us fix some positive integers $n,r$ \st   $n> r$,  a self adjoint $n\times n$ real matrix $X$ with eigenvalues 
$\la_1,\ldots, \la_n$ and some non null  real numbers   $\theta_1, \ldots, \theta_r$. We make the following hypothesis:\\ \\
(H) $\la_1, \ldots, \la_n$ are pairwise distinct  and  there are pairwise distinct indices
 $i_1, \ldots, i_{r-1}\in \{1,\ldots, n\}$ \st $\{\la_{i_1}+\theta_1, \ldots, \la_{i_{r-1}}+\theta_{r-1}\}\cap\{\la_1,\ldots, \la_n\}=\emptyset.$\\

Let us define, for $g=[g_1,\ldots, g_r]\in \mathbb{K}^{n\times r}$ , $$\widetilde{X}_g:=X+\theta_1u_1u^*_1+\cdots+\theta_ru_ru^*_r,$$ where $(u_1, \ldots, u_r)$ 
is either the orthonormalized family deduced from the columns  of $g$ by the Gram-Schmidt process or $\ff{\sqrt{n}}(g_1, \ldots, g_r)$. 

  Then the Lebesgue measure of the set of the $g$'s \st $\widetilde{X}_g$ and $X$ have at least one eigenvalue in common is null. 
\end{lem}

Now, Lemma \ref{lemdistinct} will be easy to deduce from the above.
Indeed, one can check that for $n$ large enough, $X_n$ satisfies hypothesis (H).
 We know that its eigenvalues $\la_1^{n}, \ldots, \la_n^{n}$ are distinct.
Moreover, let $\eta$ be such that  $\eta < \frac{1}{2}\min_{1\le i \le r}|\theta_i|$ and  $\eta < \frac{1}{3}\min_{i \neq j}|\ell_i-\ell_j|.$
From Assumption \ref{with}, there exists $n$ large enough so that $X_n$ has at most $p^+$ eigenvalues greater than $b+\eta,$
at most $p^-$ eigenvalues smaller than $a-\eta,$
more than $2r(p^++1)$ eigenvalues in the interval $(b-\eta, b+\eta)$
and more than $2r(p^-+1)$ eigenvalues in $(a-\eta, a+\eta).$ \\
Let us assume that $\theta_1 >0.$ Then one can find an eigenvalue $\lambda_{i_1}$ among the $p^++1$ greater ones
in   $(b-\eta, b+\eta)$ such that $\lambda_{i_1}+\theta_1$ do not belong to $\{\la_1^{n}, \ldots, \la_n^{n}\}.$
We then forget the $p^++1$ greater eigenvalues and look at the $p^++1$ following ones. Among them, one can find an eigenvalue $\lambda_{i_2}$ such that 
$\lambda_{i_2}+\theta_2$ do not belong to $\{\la_1^{n}, \ldots, \la_n^{n}\}.$ and so on.
For the negative $\theta_i$'s, we consider the  $p^-+1$ smallest eigenvalues in   $(a-\eta, a+\eta).$


We now prove Lemma \ref{vpcommon}.\\

\begin{pr} The idea of the proof is the following. We shall first prove (in Step I) that  the set of $g$'s \st $\widetilde{X}_g$ and $X$ have at least one eigenvalue in common is, up to a set of null Lebesgue measure, the set of zeroes of a polynomial function. Since it can easily be proved, by induction on the number of  variables, that the set of zeroes of any non null polynomial in several real variables 
has null Lebesgue measure, proving (in Step II) that this function is not identically null will then  imply  that the set of such $g$'s has vanishing Lebesgue measure.

Let $\beta$ be either $1$ or $2$ according to whether $\mathbb{K}$ is $\R$ or $\C$. 

{\it Step I.} Let us first treat the case where $(u_1, \ldots, u_r)=\ff{\sqrt{n}}(g_1, \ldots, g_r)$. Let us define $P$ to be the polynomial of $\beta nr$ real variables which maps $[g_1, \ldots, g_r]\in \mathbb{K}^{n\times r}$ to the resultant of the characteristic polynomials of $X$ and  $\widetilde{X}_g$. The set of $g$'s in $\mathbb{K}^{n\times r}$ \st $X$ and $\widetilde{X}_g$ have an eigenvalue in common is exactly the set of $g$'s \st $P(g)=0$ : Step I is achieved in the case where  $(u_1, \ldots, u_r)=\ff{\sqrt{n}}(g_1, \ldots, g_r)$. 

Let us now treat the case where $(u_1, \ldots, u_r)$ is  the orthonormalized family deduced from the columns  of $g$ by the Gram-Schmidt process. In this case, the resultant of the characteristic polynomials of $X$ and of $\widetilde{X}_g$ is not anymore a polynomial function of the real coordinates of $g$, so we shall use the following trick. It can easily be noticed, through  a careful look at the    Gram-Schmidt process, that 
 for all $k\in \{1,\ldots, r\}$, for all $i,j\in \{1,\ldots, n\} $, there are two   polynomial functions of  $D_k, N_{k,i,j}$ of $\beta nr$ real variables   \st the $i,j$-th entry of $u_ku_k^*$ is $\f{N_{k,i,j}(g)}{{D_{k}(g)}}$ and that $D_{k}(g)$ is positive for  any $g\in \mathbb{K}^{n\times r}$ which columns are linearly independent. Let us define   the  polynomial function of  $\beta  nr$ real variables $$D(g):=\prod_{k=1}^r D_{k}(g).$$
 For any $g$ \st $D(g)>0$ (which is the case  for  any $g\in \mathbb{K}^{n\times r}$ which columns are linearly independent), $X$ and $\widetilde{X}_g$ have no eigenvalue in common if and only if $D(g)X$ and $D(g)\widetilde{X}_g$ have no eigenvalue in common. Now, the advantage of having replaced  $X$ and $\widetilde{X}_g$  by $D(g)X$ and $D(g)\widetilde{X}_g$ is that the entries of $D(g)X$ and $D(g)\widetilde{X}_g$ are polynomial functions of $g$. Hence if one defines $P(g)$ to be the resultant of the characteristic polynomials of $D(g)X$ and $D(g)\widetilde{X}_g$, $P(g)$ is a polynomial function of the $\beta nr$ real coordinates of $g$ and, up to the set (with zero Lebesgue measure) of $g$'s in $\mathbb{K}^{n\times r}$  which columns are linearly independent, the set of $g$'s in $\mathbb{K}^{n\times r}$ \st $X$ and $\widetilde{X}_g$ have an eigenvalue in common is exactly the set of $g$'s \st $P(g)=0$ : Step I is achieved in the second case.
 
{\it Step II.}  Let us now prove that in both cases, the polynomial function $g\longmapsto P(g)$ is not identically null. To treat both cases together, it suffices to   prove that there exists $g=[g_1,\ldots, g_r]\in \mathbb{K}^{n\times r}$ with orthonormalized columns \st $\widetilde{X}_g$ and $X$ have no eigenvalue in common.  
One can suppose that $i_1=1, \ldots, i_{r-1}=r-1$, that $\la_r<\cdots<\la_n $ and that $$X=\bbm\la_1&&\\&\ddots&\\ &&\la_n\ebm.$$   We shall choose the $r-1$ first columns 
 $g_1, \ldots,g_{r-1}$ of $g$ to be the $r-1$ first elements of the canonical basis  and $g_r$ with null $r-1$ first coordinates and unit norm.  With such a choice of $g$, we have 
$$ \widetilde{X}_g=
\bbm\la_1&&&&\\
&\ddots&&&\\
&&\la_{r-1}&&\\
&&&\la_r&&\\
&&&&\ddots&\\
&&&&&\la_n
\ebm+
\bbm\theta_1&&&&&\\
&\ddots&&&&\\
&&\theta_{r-1}&&&\\
&&&&&\\
&&&&\theta_rg_rg_r^*&\\
&&&&&
\ebm.
$$ Let us suppose that $\theta_r>0$. It was shown in \cite[Section 3.2]{forrester-nagao} that
 as $g_r$  runs through the set of unit norm vectors of $\mathbb{K}^{n\times 1}$ with null $r-1$ first coordinates, 
the ordered eigenvalues  of the $n-(r-1)\times n-(r-1)$ lower right block of $ \widetilde{X}_g$ describe the set of families $\mu_r, \ldots, \mu_n$
 of real numbers which sum up to $\la_r+\cdots+\la_n+\theta_r$ and \st $$\la_r\le\mu_r\le \la_{r+1}\le\cdots\le\la_n\le\mu_n.$$
One can easily find such a family  $\mu_r, \ldots, \mu_n$ \st $$\{ \mu_r, \ldots, \mu_n\}\cap\{\la_1,\ldots, \la_n\}=\emptyset,$$
which concludes the proof, by hypothesis (H).
\end{pr}

\end{document}